\input amstex
\documentstyle{amsppt}

\catcode`\@=11
\font@\twelverm=cmr10 scaled\magstep1
\font@\twelveit=cmti10 scaled\magstep1
\font@\twelvebf=cmbx10 scaled\magstep1
\font@\twelvei=cmmi10 scaled\magstep1
\font@\twelvesy=cmsy10 scaled\magstep1
\font@\twelveex=cmex10 scaled\magstep1

\newtoks\twelvepoint@
\def\twelvepoint{\normalbaselineskip15\p@
 \abovedisplayskip15\p@ plus3.6\p@ minus10.8\p@
 \belowdisplayskip\abovedisplayskip
 \abovedisplayshortskip\z@ plus3.6\p@
 \belowdisplayshortskip8.4\p@ plus3.6\p@ minus4.8\p@
 \textonlyfont@\rm\twelverm \textonlyfont@\it\twelveit
 \textonlyfont@\sl\twelvesl \textonlyfont@\bf\twelvebf
 \textonlyfont@\smc\twelvesmc \textonlyfont@\tt\twelvett
%
 \ifsyntax@ \def\big##1{{\hbox{$\left##1\right.$}}}%
  \let\Big\big \let\bigg\big \let\Bigg\big
 \else
  \textfont\z@=\twelverm  \scriptfont\z@=\tenrm  \scriptscriptfont\z@=\sevenrm
  \textfont\@ne=\twelvei  \scriptfont\@ne=\teni  \scriptscriptfont\@ne=\seveni
  \textfont\tw@=\twelvesy \scriptfont\tw@=\tensy \scriptscriptfont\tw@=\sevensy
  \textfont\thr@@=\twelveex \scriptfont\thr@@=\tenex
        \scriptscriptfont\thr@@=\tenex
  \textfont\itfam=\twelveit \scriptfont\itfam=\tenit
        \scriptscriptfont\itfam=\tenit
  \textfont\bffam=\twelvebf \scriptfont\bffam=\tenbf
        \scriptscriptfont\bffam=\sevenbf
  \setbox\strutbox\hbox{\vrule height10.2\p@ depth4.2\p@ width\z@}%
  \setbox\strutbox@\hbox{\lower.6\normallineskiplimit\vbox{%
        \kern-\normallineskiplimit\copy\strutbox}}%
 \setbox\z@\vbox{\hbox{$($}\kern\z@}\bigsize@=1.4\ht\z@
 \fi
 \normalbaselines\rm\ex@.2326ex\jot3.6\ex@\the\twelvepoint@}

\font@\fourteenrm=cmr10 scaled\magstep2
\font@\fourteenit=cmti10 scaled\magstep2
\font@\fourteensl=cmsl10 scaled\magstep2
\font@\fourteensmc=cmcsc10 scaled\magstep2
\font@\fourteentt=cmtt10 scaled\magstep2
\font@\fourteenbf=cmbx10 scaled\magstep2
\font@\fourteeni=cmmi10 scaled\magstep2
\font@\fourteensy=cmsy10 scaled\magstep2
\font@\fourteenex=cmex10 scaled\magstep2
\font@\fourteenmsa=msam10 scaled\magstep2
\font@\fourteeneufm=eufm10 scaled\magstep2
\font@\fourteenmsb=msbm10 scaled\magstep2
\newtoks\fourteenpoint@
\def\fourteenpoint{\normalbaselineskip15\p@
 \abovedisplayskip18\p@ plus4.3\p@ minus12.9\p@
 \belowdisplayskip\abovedisplayskip
 \abovedisplayshortskip\z@ plus4.3\p@
 \belowdisplayshortskip10.1\p@ plus4.3\p@ minus5.8\p@
 \textonlyfont@\rm\fourteenrm \textonlyfont@\it\fourteenit
 \textonlyfont@\sl\fourteensl \textonlyfont@\bf\fourteenbf
 \textonlyfont@\smc\fourteensmc \textonlyfont@\tt\fourteentt
%
 \ifsyntax@ \def\big##1{{\hbox{$\left##1\right.$}}}%
  \let\Big\big \let\bigg\big \let\Bigg\big
 \else
  \textfont\z@=\fourteenrm  \scriptfont\z@=\twelverm  \scriptscriptfont\z@=\tenrm
  \textfont\@ne=\fourteeni  \scriptfont\@ne=\twelvei  \scriptscriptfont\@ne=\teni
  \textfont\tw@=\fourteensy \scriptfont\tw@=\twelvesy \scriptscriptfont\tw@=\tensy
  \textfont\thr@@=\fourteenex \scriptfont\thr@@=\twelveex
        \scriptscriptfont\thr@@=\twelveex
  \textfont\itfam=\fourteenit \scriptfont\itfam=\twelveit
        \scriptscriptfont\itfam=\twelveit
  \textfont\bffam=\fourteenbf \scriptfont\bffam=\twelvebf
        \scriptscriptfont\bffam=\tenbf
  \setbox\strutbox\hbox{\vrule height12.2\p@ depth5\p@ width\z@}%
  \setbox\strutbox@\hbox{\lower.72\normallineskiplimit\vbox{%
        \kern-\normallineskiplimit\copy\strutbox}}%
 \setbox\z@\vbox{\hbox{$($}\kern\z@}\bigsize@=1.7\ht\z@
 \fi
 \normalbaselines\rm\ex@.2326ex\jot4.3\ex@\the\fourteenpoint@}

\font@\seventeenrm=cmr10 scaled\magstep3
\font@\seventeenit=cmti10 scaled\magstep3
\font@\seventeensl=cmsl10 scaled\magstep3
\font@\seventeensmc=cmcsc10 scaled\magstep3
\font@\seventeentt=cmtt10 scaled\magstep3
\font@\seventeenbf=cmbx10 scaled\magstep3
\font@\seventeeni=cmmi10 scaled\magstep3
\font@\seventeensy=cmsy10 scaled\magstep3
\font@\seventeenex=cmex10 scaled\magstep3
\font@\seventeenmsa=msam10 scaled\magstep3
\font@\seventeeneufm=eufm10 scaled\magstep3
\font@\seventeenmsb=msbm10 scaled\magstep3
\newtoks\seventeenpoint@
\def\seventeenpoint{\normalbaselineskip18\p@
 \abovedisplayskip21.6\p@ plus5.2\p@ minus15.4\p@
 \belowdisplayskip\abovedisplayskip
 \abovedisplayshortskip\z@ plus5.2\p@
 \belowdisplayshortskip12.1\p@ plus5.2\p@ minus7\p@
 \textonlyfont@\rm\seventeenrm \textonlyfont@\it\seventeenit
 \textonlyfont@\sl\seventeensl \textonlyfont@\bf\seventeenbf
 \textonlyfont@\smc\seventeensmc \textonlyfont@\tt\seventeentt
%
 \ifsyntax@ \def\big##1{{\hbox{$\left##1\right.$}}}%
  \let\Big\big \let\bigg\big \let\Bigg\big
 \else
  \textfont\z@=\seventeenrm  \scriptfont\z@=\fourteenrm  \scriptscriptfont\z@=\twelverm
  \textfont\@ne=\seventeeni  \scriptfont\@ne=\fourteeni  \scriptscriptfont\@ne=\twelvei
  \textfont\tw@=\seventeensy \scriptfont\tw@=\fourteensy \scriptscriptfont\tw@=\twelvesy
  \textfont\thr@@=\seventeenex \scriptfont\thr@@=\fourteenex
        \scriptscriptfont\thr@@=\fourteenex
  \textfont\itfam=\seventeenit \scriptfont\itfam=\fourteenit
        \scriptscriptfont\itfam=\fourteenit
  \textfont\bffam=\seventeenbf \scriptfont\bffam=\fourteenbf
        \scriptscriptfont\bffam=\twelvebf
  \setbox\strutbox\hbox{\vrule height14.6\p@ depth6\p@ width\z@}%
  \setbox\strutbox@\hbox{\lower.86\normallineskiplimit\vbox{%
        \kern-\normallineskiplimit\copy\strutbox}}%
 \setbox\z@\vbox{\hbox{$($}\kern\z@}\bigsize@=2\ht\z@
 \fi
 \normalbaselines\rm\ex@.2326ex\jot5.2\ex@\the\seventeenpoint@}

\catcode`\@=13
\catcode`\@=11
\font\tenln    = line10
\font\tenlnw   = linew10

\newskip\Einheit \Einheit=0.5cm
\newcount\xcoord \newcount\ycoord
\newdimen\xdim \newdimen\ydim \newdimen\PfadD@cke \newdimen\Pfadd@cke

\newcount\@tempcnta
\newcount\@tempcntb

\newdimen\@tempdima
\newdimen\@tempdimb

\newdimen\@wholewidth
\newdimen\@halfwidth

\newcount\@xarg
\newcount\@yarg
\newcount\@yyarg
\newbox\@linechar
\newbox\@tempboxa
\newdimen\@linelen
\newdimen\@clnwd
\newdimen\@clnht

\newif\if@negarg

\def\@whilenoop#1{}
\def\@whiledim#1\do #2{\ifdim #1\relax#2\@iwhiledim{#1\relax#2}\fi}
\def\@iwhiledim#1{\ifdim #1\let\@nextwhile=\@iwhiledim
        \else\let\@nextwhile=\@whilenoop\fi\@nextwhile{#1}}

\def\@whileswnoop#1\fi{}
\def\@whilesw#1\fi#2{#1#2\@iwhilesw{#1#2}\fi\fi}
\def\@iwhilesw#1\fi{#1\let\@nextwhile=\@iwhilesw
         \else\let\@nextwhile=\@whileswnoop\fi\@nextwhile{#1}\fi}

\def\thinlines{\let\@linefnt\tenln \let\@circlefnt\tencirc
  \@wholewidth\fontdimen8\tenln \@halfwidth .5\@wholewidth}
\def\thicklines{\let\@linefnt\tenlnw \let\@circlefnt\tencircw
  \@wholewidth\fontdimen8\tenlnw \@halfwidth .5\@wholewidth}
\thinlines

\PfadD@cke1pt \Pfadd@cke0.5pt
\def\PfadDicke#1{\PfadD@cke#1 \divide\PfadD@cke by2 \Pfadd@cke\PfadD@cke \multiply\PfadD@cke by2}
\long\def\LOOP#1\REPEAT{\def\BODY{#1}\ITERATE}
\def\ITERATE{\BODY \let\next\ITERATE \else\let\next\relax\fi \next}
\let\REPEAT=\fi
\def\Punkt{\hbox{\raise-2pt\hbox to0pt{\hss$\ssize\bullet$\hss}}}
\def\DuennPunkt(#1,#2){\unskip
  \raise#2 \Einheit\hbox to0pt{\hskip#1 \Einheit
          \raise-2.5pt\hbox to0pt{\hss$\bullet$\hss}\hss}}
\def\NormalPunkt(#1,#2){\unskip
  \raise#2 \Einheit\hbox to0pt{\hskip#1 \Einheit
          \raise-3pt\hbox to0pt{\hss\twelvepoint$\bullet$\hss}\hss}}
\def\DickPunkt(#1,#2){\unskip
  \raise#2 \Einheit\hbox to0pt{\hskip#1 \Einheit
          \raise-4pt\hbox to0pt{\hss\fourteenpoint$\bullet$\hss}\hss}}
\def\Kreis(#1,#2){\unskip
  \raise#2 \Einheit\hbox to0pt{\hskip#1 \Einheit
          \raise-4pt\hbox to0pt{\hss\fourteenpoint$\circ$\hss}\hss}}

\def\Line@(#1,#2)#3{\@xarg #1\relax \@yarg #2\relax
\@linelen=#3\Einheit
\ifnum\@xarg =0 \@vline
  \else \ifnum\@yarg =0 \@hline \else \@sline\fi
\fi}

\def\@sline{\ifnum\@xarg< 0 \@negargtrue \@xarg -\@xarg \@yyarg -\@yarg
  \else \@negargfalse \@yyarg \@yarg \fi
\ifnum \@yyarg >0 \@tempcnta\@yyarg \else \@tempcnta -\@yyarg \fi
\ifnum\@tempcnta>6 \@badlinearg\@tempcnta0 \fi
\ifnum\@xarg>6 \@badlinearg\@xarg 1 \fi
\setbox\@linechar\hbox{\@linefnt\@getlinechar(\@xarg,\@yyarg)}%
\ifnum \@yarg >0 \let\@upordown\raise \@clnht\z@
   \else\let\@upordown\lower \@clnht \ht\@linechar\fi
\@clnwd=\wd\@linechar
\if@negarg \hskip -\wd\@linechar \def\@tempa{\hskip -2\wd\@linechar}\else
     \let\@tempa\relax \fi
\@whiledim \@clnwd <\@linelen \do
  {\@upordown\@clnht\copy\@linechar
   \@tempa
   \advance\@clnht \ht\@linechar
   \advance\@clnwd \wd\@linechar}%
\advance\@clnht -\ht\@linechar
\advance\@clnwd -\wd\@linechar
\@tempdima\@linelen\advance\@tempdima -\@clnwd
\@tempdimb\@tempdima\advance\@tempdimb -\wd\@linechar
\if@negarg \hskip -\@tempdimb \else \hskip \@tempdimb \fi
\multiply\@tempdima \@m
\@tempcnta \@tempdima \@tempdima \wd\@linechar \divide\@tempcnta \@tempdima
\@tempdima \ht\@linechar \multiply\@tempdima \@tempcnta
\divide\@tempdima \@m
\advance\@clnht \@tempdima
\ifdim \@linelen <\wd\@linechar
   \hskip \wd\@linechar
  \else\@upordown\@clnht\copy\@linechar\fi}

\def\@hline{\ifnum \@xarg <0 \hskip -\@linelen \fi
\vrule height\Pfadd@cke width \@linelen depth\Pfadd@cke
\ifnum \@xarg <0 \hskip -\@linelen \fi}

\def\@getlinechar(#1,#2){\@tempcnta#1\relax\multiply\@tempcnta 8
\advance\@tempcnta -9 \ifnum #2>0 \advance\@tempcnta #2\relax\else
\advance\@tempcnta -#2\relax\advance\@tempcnta 64 \fi
\char\@tempcnta}

\def\Vektor(#1,#2)#3(#4,#5){\unskip\leavevmode
  \xcoord#4\relax \ycoord#5\relax
      \raise\ycoord \Einheit\hbox to0pt{\hskip\xcoord \Einheit
         \Vector@(#1,#2){#3}\hss}}

\def\Vector@(#1,#2)#3{\@xarg #1\relax \@yarg #2\relax
\@tempcnta \ifnum\@xarg<0 -\@xarg\else\@xarg\fi
\ifnum\@tempcnta<5\relax
\@linelen=#3\Einheit
\ifnum\@xarg =0 \@vvector
  \else \ifnum\@yarg =0 \@hvector \else \@svector\fi
\fi
\else\@badlinearg\fi}

\def\@hvector{\@hline\hbox to 0pt{\@linefnt
\ifnum \@xarg <0 \@getlarrow(1,0)\hss\else
    \hss\@getrarrow(1,0)\fi}}

\def\@vvector{\ifnum \@yarg <0 \@downvector \else \@upvector \fi}

\def\@svector{\@sline
\@tempcnta\@yarg \ifnum\@tempcnta <0 \@tempcnta=-\@tempcnta\fi
\ifnum\@tempcnta <5
  \hskip -\wd\@linechar
  \@upordown\@clnht \hbox{\@linefnt  \if@negarg
  \@getlarrow(\@xarg,\@yyarg) \else \@getrarrow(\@xarg,\@yyarg) \fi}%
\else\@badlinearg\fi}

\def\@upline{\hbox to \z@{\hskip -.5\Pfadd@cke \vrule width \Pfadd@cke
   height \@linelen depth \z@\hss}}

\def\@downline{\hbox to \z@{\hskip -.5\Pfadd@cke \vrule width \Pfadd@cke
   height \z@ depth \@linelen \hss}}

\def\@upvector{\@upline\setbox\@tempboxa\hbox{\@linefnt\char'66}\raise
     \@linelen \hbox to\z@{\lower \ht\@tempboxa\box\@tempboxa\hss}}

\def\@downvector{\@downline\lower \@linelen
      \hbox to \z@{\@linefnt\char'77\hss}}

\def\@getlarrow(#1,#2){\ifnum #2 =\z@ \@tempcnta='33\else
\@tempcnta=#1\relax\multiply\@tempcnta \sixt@@n \advance\@tempcnta
-9 \@tempcntb=#2\relax\multiply\@tempcntb \tw@
\ifnum \@tempcntb >0 \advance\@tempcnta \@tempcntb\relax
\else\advance\@tempcnta -\@tempcntb\advance\@tempcnta 64
\fi\fi\char\@tempcnta}

\def\@getrarrow(#1,#2){\@tempcntb=#2\relax
\ifnum\@tempcntb < 0 \@tempcntb=-\@tempcntb\relax\fi
\ifcase \@tempcntb\relax \@tempcnta='55 \or
\ifnum #1<3 \@tempcnta=#1\relax\multiply\@tempcnta
24 \advance\@tempcnta -6 \else \ifnum #1=3 \@tempcnta=49
\else\@tempcnta=58 \fi\fi\or
\ifnum #1<3 \@tempcnta=#1\relax\multiply\@tempcnta
24 \advance\@tempcnta -3 \else \@tempcnta=51\fi\or
\@tempcnta=#1\relax\multiply\@tempcnta
\sixt@@n \advance\@tempcnta -\tw@ \else
\@tempcnta=#1\relax\multiply\@tempcnta
\sixt@@n \advance\@tempcnta 7 \fi\ifnum #2<0 \advance\@tempcnta 64 \fi
\char\@tempcnta}

\def\Diagonale(#1,#2)#3{\unskip\leavevmode
  \xcoord#1\relax \ycoord#2\relax
      \raise\ycoord \Einheit\hbox to0pt{\hskip\xcoord \Einheit
         \Line@(1,1){#3}\hss}}
\def\AntiDiagonale(#1,#2)#3{\unskip\leavevmode
  \xcoord#1\relax \ycoord#2\relax 
      \raise\ycoord \Einheit\hbox to0pt{\hskip\xcoord \Einheit
         \Line@(1,-1){#3}\hss}}
\def\Pfad(#1,#2),#3\endPfad{\unskip\leavevmode
  \xcoord#1 \ycoord#2 \thicklines\ZeichnePfad#3\endPfad\thinlines}
\def\ZeichnePfad#1{\ifx#1\endPfad\let\next\relax
  \else\let\next\ZeichnePfad
    \ifnum#1=1
      \raise\ycoord \Einheit\hbox to0pt{\hskip\xcoord \Einheit
         \vrule height\Pfadd@cke width1 \Einheit depth\Pfadd@cke\hss}%
      \advance\xcoord by 1
    \else\ifnum#1=2
      \raise\ycoord \Einheit\hbox to0pt{\hskip\xcoord \Einheit
        \hbox{\hskip-\PfadD@cke\vrule height1 \Einheit width\PfadD@cke depth0pt}\hss}%
      \advance\ycoord by 1
    \else\ifnum#1=3
      \raise\ycoord \Einheit\hbox to0pt{\hskip\xcoord \Einheit
         \Line@(1,1){1}\hss}
      \advance\xcoord by 1
      \advance\ycoord by 1
    \else\ifnum#1=4
      \raise\ycoord \Einheit\hbox to0pt{\hskip\xcoord \Einheit
         \Line@(1,-1){1}\hss}
      \advance\xcoord by 1
      \advance\ycoord by -1
    \fi\fi\fi\fi
  \fi\next}
\def\hSSchritt{\leavevmode\raise-.4pt\hbox to0pt{\hss.\hss}\hskip.2\Einheit
  \raise-.4pt\hbox to0pt{\hss.\hss}\hskip.2\Einheit
  \raise-.4pt\hbox to0pt{\hss.\hss}\hskip.2\Einheit
  \raise-.4pt\hbox to0pt{\hss.\hss}\hskip.2\Einheit
  \raise-.4pt\hbox to0pt{\hss.\hss}\hskip.2\Einheit}
\def\vSSchritt{\vbox{\baselineskip.2\Einheit\lineskiplimit0pt
\hbox{.}\hbox{.}\hbox{.}\hbox{.}\hbox{.}}}
\def\DSSchritt{\leavevmode\raise-.4pt\hbox to0pt{%
  \hbox to0pt{\hss.\hss}\hskip.2\Einheit
  \raise.2\Einheit\hbox to0pt{\hss.\hss}\hskip.2\Einheit
  \raise.4\Einheit\hbox to0pt{\hss.\hss}\hskip.2\Einheit
  \raise.6\Einheit\hbox to0pt{\hss.\hss}\hskip.2\Einheit
  \raise.8\Einheit\hbox to0pt{\hss.\hss}\hss}}
\def\dSSchritt{\leavevmode\raise-.4pt\hbox to0pt{%
  \hbox to0pt{\hss.\hss}\hskip.2\Einheit
  \raise-.2\Einheit\hbox to0pt{\hss.\hss}\hskip.2\Einheit
  \raise-.4\Einheit\hbox to0pt{\hss.\hss}\hskip.2\Einheit
  \raise-.6\Einheit\hbox to0pt{\hss.\hss}\hskip.2\Einheit
  \raise-.8\Einheit\hbox to0pt{\hss.\hss}\hss}}
\def\SPfad(#1,#2),#3\endSPfad{\unskip\leavevmode
  \xcoord#1 \ycoord#2 \ZeichneSPfad#3\endSPfad}
\def\ZeichneSPfad#1{\ifx#1\endSPfad\let\next\relax
  \else\let\next\ZeichneSPfad
    \ifnum#1=1
      \raise\ycoord \Einheit\hbox to0pt{\hskip\xcoord \Einheit
         \hSSchritt\hss}%
      \advance\xcoord by 1
    \else\ifnum#1=2
      \raise\ycoord \Einheit\hbox to0pt{\hskip\xcoord \Einheit
        \hbox{\hskip-2pt \vSSchritt}\hss}%
      \advance\ycoord by 1
    \else\ifnum#1=3
      \raise\ycoord \Einheit\hbox to0pt{\hskip\xcoord \Einheit
         \DSSchritt\hss}
      \advance\xcoord by 1
      \advance\ycoord by 1
    \else\ifnum#1=4
      \raise\ycoord \Einheit\hbox to0pt{\hskip\xcoord \Einheit
         \dSSchritt\hss}
      \advance\xcoord by 1
      \advance\ycoord by -1
    \fi\fi\fi\fi
  \fi\next}
\def\Koordinatenachsen(#1,#2){\unskip
 \hbox to0pt{\hskip-.5pt\vrule height#2 \Einheit width.5pt depth1 \Einheit}%
 \hbox to0pt{\hskip-1 \Einheit \xcoord#1 \advance\xcoord by1
    \vrule height0.25pt width\xcoord \Einheit depth0.25pt\hss}}
\def\Koordinatenachsen(#1,#2)(#3,#4){\unskip
 \hbox to0pt{\hskip-.5pt \ycoord-#4 \advance\ycoord by1
    \vrule height#2 \Einheit width.5pt depth\ycoord \Einheit}%
 \hbox to0pt{\hskip-1 \Einheit \hskip#3\Einheit 
    \xcoord#1 \advance\xcoord by1 \advance\xcoord by-#3 
    \vrule height0.25pt width\xcoord \Einheit depth0.25pt\hss}}
\def\Gitter(#1,#2){\unskip \xcoord0 \ycoord0 \leavevmode
  \LOOP\ifnum\ycoord<#2
    \loop\ifnum\xcoord<#1
      \raise\ycoord \Einheit\hbox to0pt{\hskip\xcoord \Einheit\Punkt\hss}%
      \advance\xcoord by1
    \repeat
    \xcoord0
    \advance\ycoord by1
  \REPEAT}
\def\Gitter(#1,#2)(#3,#4){\unskip \xcoord#3 \ycoord#4 \leavevmode
  \LOOP\ifnum\ycoord<#2
    \loop\ifnum\xcoord<#1
      \raise\ycoord \Einheit\hbox to0pt{\hskip\xcoord \Einheit\Punkt\hss}%
      \advance\xcoord by1
    \repeat
    \xcoord#3
    \advance\ycoord by1
  \REPEAT}
\def\Label#1#2(#3,#4){\unskip \xdim#3 \Einheit \ydim#4 \Einheit
  \def\lo{\advance\xdim by-.5 \Einheit \advance\ydim by.5 \Einheit}%
  \def\llo{\advance\xdim by-.25cm \advance\ydim by.5 \Einheit}%
  \def\loo{\advance\xdim by-.5 \Einheit \advance\ydim by.25cm}%
  \def\o{\advance\ydim by.25cm}%
  \def\ro{\advance\xdim by.5 \Einheit \advance\ydim by.5 \Einheit}%
  \def\rro{\advance\xdim by.25cm \advance\ydim by.5 \Einheit}%
  \def\roo{\advance\xdim by.5 \Einheit \advance\ydim by.25cm}%
  \def\l{\advance\xdim by-.30cm}%
  \def\r{\advance\xdim by.30cm}%
  \def\lu{\advance\xdim by-.5 \Einheit \advance\ydim by-.6 \Einheit}%
  \def\llu{\advance\xdim by-.25cm \advance\ydim by-.6 \Einheit}%
  \def\luu{\advance\xdim by-.5 \Einheit \advance\ydim by-.30cm}%
  \def\u{\advance\ydim by-.30cm}%
  \def\ru{\advance\xdim by.5 \Einheit \advance\ydim by-.6 \Einheit}%
  \def\rru{\advance\xdim by.25cm \advance\ydim by-.6 \Einheit}%
  \def\ruu{\advance\xdim by.5 \Einheit \advance\ydim by-.30cm}%
  #1\raise\ydim\hbox to0pt{\hskip\xdim
     \vbox to0pt{\vss\hbox to0pt{\hss$#2$\hss}\vss}\hss}%
}
\catcode`\@=13

\catcode`\@=11
\font\tenln    = line10
\font\tenlnw   = linew10

\newskip\Einheit \Einheit=0.5cm
\newcount\xcoord \newcount\ycoord
\newdimen\xdim \newdimen\ydim \newdimen\PfadD@cke \newdimen\Pfadd@cke

\newcount\@tempcnta
\newcount\@tempcntb

\newdimen\@tempdima
\newdimen\@tempdimb

\newdimen\@wholewidth
\newdimen\@halfwidth

\newcount\@xarg
\newcount\@yarg
\newcount\@yyarg
\newbox\@linechar
\newbox\@tempboxa
\newdimen\@linelen
\newdimen\@clnwd
\newdimen\@clnht

\newif\if@negarg

\def\@whilenoop#1{}
\def\@whiledim#1\do #2{\ifdim #1\relax#2\@iwhiledim{#1\relax#2}\fi}
\def\@iwhiledim#1{\ifdim #1\let\@nextwhile=\@iwhiledim
        \else\let\@nextwhile=\@whilenoop\fi\@nextwhile{#1}}

\def\@whileswnoop#1\fi{}
\def\@whilesw#1\fi#2{#1#2\@iwhilesw{#1#2}\fi\fi}
\def\@iwhilesw#1\fi{#1\let\@nextwhile=\@iwhilesw
         \else\let\@nextwhile=\@whileswnoop\fi\@nextwhile{#1}\fi}

\def\thinlines{\let\@linefnt\tenln \let\@circlefnt\tencirc
  \@wholewidth\fontdimen8\tenln \@halfwidth .5\@wholewidth}
\def\thicklines{\let\@linefnt\tenlnw \let\@circlefnt\tencircw
  \@wholewidth\fontdimen8\tenlnw \@halfwidth .5\@wholewidth}
\thinlines

\PfadD@cke1pt \Pfadd@cke0.5pt
\def\PfadDicke#1{\PfadD@cke#1 \divide\PfadD@cke by2 \Pfadd@cke\PfadD@cke \multiply\PfadD@cke by2}
\long\def\LOOP#1\REPEAT{\def\BODY{#1}\ITERATE}
\def\ITERATE{\BODY \let\next\ITERATE \else\let\next\relax\fi \next}
\let\REPEAT=\fi
\def\Punkt{\hbox{\raise-2pt\hbox to0pt{\hss$\ssize\bullet$\hss}}}
\def\DuennPunkt(#1,#2){\unskip
  \raise#2 \Einheit\hbox to0pt{\hskip#1 \Einheit
          \raise-2.5pt\hbox to0pt{\hss$\bullet$\hss}\hss}}
\def\NormalPunkt(#1,#2){\unskip
  \raise#2 \Einheit\hbox to0pt{\hskip#1 \Einheit
          \raise-3pt\hbox to0pt{\hss\twelvepoint$\bullet$\hss}\hss}}
\def\DickPunkt(#1,#2){\unskip
  \raise#2 \Einheit\hbox to0pt{\hskip#1 \Einheit
          \raise-4pt\hbox to0pt{\hss\fourteenpoint$\bullet$\hss}\hss}}
\def\Kreis(#1,#2){\unskip
  \raise#2 \Einheit\hbox to0pt{\hskip#1 \Einheit
          \raise-4pt\hbox to0pt{\hss\fourteenpoint$\circ$\hss}\hss}}

\def\Line@(#1,#2)#3{\@xarg #1\relax \@yarg #2\relax
\@linelen=#3\Einheit
\ifnum\@xarg =0 \@vline
  \else \ifnum\@yarg =0 \@hline \else \@sline\fi
\fi}

\def\@sline{\ifnum\@xarg< 0 \@negargtrue \@xarg -\@xarg \@yyarg -\@yarg
  \else \@negargfalse \@yyarg \@yarg \fi
\ifnum \@yyarg >0 \@tempcnta\@yyarg \else \@tempcnta -\@yyarg \fi
\ifnum\@tempcnta>6 \@badlinearg\@tempcnta0 \fi
\ifnum\@xarg>6 \@badlinearg\@xarg 1 \fi
\setbox\@linechar\hbox{\@linefnt\@getlinechar(\@xarg,\@yyarg)}%
\ifnum \@yarg >0 \let\@upordown\raise \@clnht\z@
   \else\let\@upordown\lower \@clnht \ht\@linechar\fi
\@clnwd=\wd\@linechar
\if@negarg \hskip -\wd\@linechar \def\@tempa{\hskip -2\wd\@linechar}\else
     \let\@tempa\relax \fi
\@whiledim \@clnwd <\@linelen \do
  {\@upordown\@clnht\copy\@linechar
   \@tempa
   \advance\@clnht \ht\@linechar
   \advance\@clnwd \wd\@linechar}%
\advance\@clnht -\ht\@linechar
\advance\@clnwd -\wd\@linechar
\@tempdima\@linelen\advance\@tempdima -\@clnwd
\@tempdimb\@tempdima\advance\@tempdimb -\wd\@linechar
\if@negarg \hskip -\@tempdimb \else \hskip \@tempdimb \fi
\multiply\@tempdima \@m
\@tempcnta \@tempdima \@tempdima \wd\@linechar \divide\@tempcnta \@tempdima
\@tempdima \ht\@linechar \multiply\@tempdima \@tempcnta
\divide\@tempdima \@m
\advance\@clnht \@tempdima
\ifdim \@linelen <\wd\@linechar
   \hskip \wd\@linechar
  \else\@upordown\@clnht\copy\@linechar\fi}

\def\@hline{\ifnum \@xarg <0 \hskip -\@linelen \fi
\vrule height\Pfadd@cke width \@linelen depth\Pfadd@cke
\ifnum \@xarg <0 \hskip -\@linelen \fi}

\def\@getlinechar(#1,#2){\@tempcnta#1\relax\multiply\@tempcnta 8
\advance\@tempcnta -9 \ifnum #2>0 \advance\@tempcnta #2\relax\else
\advance\@tempcnta -#2\relax\advance\@tempcnta 64 \fi
\char\@tempcnta}

\def\Vektor(#1,#2)#3(#4,#5){\unskip\leavevmode
  \xcoord#4\relax \ycoord#5\relax
      \raise\ycoord \Einheit\hbox to0pt{\hskip\xcoord \Einheit
         \Vector@(#1,#2){#3}\hss}}

\def\Vector@(#1,#2)#3{\@xarg #1\relax \@yarg #2\relax
\@tempcnta \ifnum\@xarg<0 -\@xarg\else\@xarg\fi
\ifnum\@tempcnta<5\relax
\@linelen=#3\Einheit
\ifnum\@xarg =0 \@vvector
  \else \ifnum\@yarg =0 \@hvector \else \@svector\fi
\fi
\else\@badlinearg\fi}

\def\@hvector{\@hline\hbox to 0pt{\@linefnt
\ifnum \@xarg <0 \@getlarrow(1,0)\hss\else
    \hss\@getrarrow(1,0)\fi}}

\def\@vvector{\ifnum \@yarg <0 \@downvector \else \@upvector \fi}

\def\@svector{\@sline
\@tempcnta\@yarg \ifnum\@tempcnta <0 \@tempcnta=-\@tempcnta\fi
\ifnum\@tempcnta <5
  \hskip -\wd\@linechar
  \@upordown\@clnht \hbox{\@linefnt  \if@negarg
  \@getlarrow(\@xarg,\@yyarg) \else \@getrarrow(\@xarg,\@yyarg) \fi}%
\else\@badlinearg\fi}

\def\@upline{\hbox to \z@{\hskip -.5\Pfadd@cke \vrule width \Pfadd@cke
   height \@linelen depth \z@\hss}}

\def\@downline{\hbox to \z@{\hskip -.5\Pfadd@cke \vrule width \Pfadd@cke
   height \z@ depth \@linelen \hss}}

\def\@upvector{\@upline\setbox\@tempboxa\hbox{\@linefnt\char'66}\raise
     \@linelen \hbox to\z@{\lower \ht\@tempboxa\box\@tempboxa\hss}}

\def\@downvector{\@downline\lower \@linelen
      \hbox to \z@{\@linefnt\char'77\hss}}

\def\@getlarrow(#1,#2){\ifnum #2 =\z@ \@tempcnta='33\else
\@tempcnta=#1\relax\multiply\@tempcnta \sixt@@n \advance\@tempcnta
-9 \@tempcntb=#2\relax\multiply\@tempcntb \tw@
\ifnum \@tempcntb >0 \advance\@tempcnta \@tempcntb\relax
\else\advance\@tempcnta -\@tempcntb\advance\@tempcnta 64
\fi\fi\char\@tempcnta}

\def\@getrarrow(#1,#2){\@tempcntb=#2\relax
\ifnum\@tempcntb < 0 \@tempcntb=-\@tempcntb\relax\fi
\ifcase \@tempcntb\relax \@tempcnta='55 \or
\ifnum #1<3 \@tempcnta=#1\relax\multiply\@tempcnta
24 \advance\@tempcnta -6 \else \ifnum #1=3 \@tempcnta=49
\else\@tempcnta=58 \fi\fi\or
\ifnum #1<3 \@tempcnta=#1\relax\multiply\@tempcnta
24 \advance\@tempcnta -3 \else \@tempcnta=51\fi\or
\@tempcnta=#1\relax\multiply\@tempcnta
\sixt@@n \advance\@tempcnta -\tw@ \else
\@tempcnta=#1\relax\multiply\@tempcnta
\sixt@@n \advance\@tempcnta 7 \fi\ifnum #2<0 \advance\@tempcnta 64 \fi
\char\@tempcnta}

\def\Diagonale(#1,#2)#3{\unskip\leavevmode
  \xcoord#1\relax \ycoord#2\relax
      \raise\ycoord \Einheit\hbox to0pt{\hskip\xcoord \Einheit
         \Line@(1,1){#3}\hss}}
\def\AntiDiagonale(#1,#2)#3{\unskip\leavevmode
  \xcoord#1\relax \ycoord#2\relax 
      \raise\ycoord \Einheit\hbox to0pt{\hskip\xcoord \Einheit
         \Line@(1,-1){#3}\hss}}
\def\Pfad(#1,#2),#3\endPfad{\unskip\leavevmode
  \xcoord#1 \ycoord#2 \thicklines\ZeichnePfad#3\endPfad\thinlines}
\def\ZeichnePfad#1{\ifx#1\endPfad\let\next\relax
  \else\let\next\ZeichnePfad
    \ifnum#1=1
      \raise\ycoord \Einheit\hbox to0pt{\hskip\xcoord \Einheit
         \vrule height\Pfadd@cke width1 \Einheit depth\Pfadd@cke\hss}%
      \advance\xcoord by 1
    \else\ifnum#1=2
      \raise\ycoord \Einheit\hbox to0pt{\hskip\xcoord \Einheit
        \hbox{\hskip-\PfadD@cke\vrule height1 \Einheit width\PfadD@cke depth0pt}\hss}%
      \advance\ycoord by 1
    \else\ifnum#1=3
      \raise\ycoord \Einheit\hbox to0pt{\hskip\xcoord \Einheit
         \Line@(1,1){1}\hss}
      \advance\xcoord by 1
      \advance\ycoord by 1
    \else\ifnum#1=4
      \raise\ycoord \Einheit\hbox to0pt{\hskip\xcoord \Einheit
         \Line@(1,-1){1}\hss}
      \advance\xcoord by 1
      \advance\ycoord by -1
    \fi\fi\fi\fi
  \fi\next}
\def\hSSchritt{\leavevmode\raise-.4pt\hbox to0pt{\hss.\hss}\hskip.2\Einheit
  \raise-.4pt\hbox to0pt{\hss.\hss}\hskip.2\Einheit
  \raise-.4pt\hbox to0pt{\hss.\hss}\hskip.2\Einheit
  \raise-.4pt\hbox to0pt{\hss.\hss}\hskip.2\Einheit
  \raise-.4pt\hbox to0pt{\hss.\hss}\hskip.2\Einheit}
\def\vSSchritt{\vbox{\baselineskip.2\Einheit\lineskiplimit0pt
\hbox{.}\hbox{.}\hbox{.}\hbox{.}\hbox{.}}}
\def\DSSchritt{\leavevmode\raise-.4pt\hbox to0pt{%
  \hbox to0pt{\hss.\hss}\hskip.2\Einheit
  \raise.2\Einheit\hbox to0pt{\hss.\hss}\hskip.2\Einheit
  \raise.4\Einheit\hbox to0pt{\hss.\hss}\hskip.2\Einheit
  \raise.6\Einheit\hbox to0pt{\hss.\hss}\hskip.2\Einheit
  \raise.8\Einheit\hbox to0pt{\hss.\hss}\hss}}
\def\dSSchritt{\leavevmode\raise-.4pt\hbox to0pt{%
  \hbox to0pt{\hss.\hss}\hskip.2\Einheit
  \raise-.2\Einheit\hbox to0pt{\hss.\hss}\hskip.2\Einheit
  \raise-.4\Einheit\hbox to0pt{\hss.\hss}\hskip.2\Einheit
  \raise-.6\Einheit\hbox to0pt{\hss.\hss}\hskip.2\Einheit
  \raise-.8\Einheit\hbox to0pt{\hss.\hss}\hss}}
\def\SPfad(#1,#2),#3\endSPfad{\unskip\leavevmode
  \xcoord#1 \ycoord#2 \ZeichneSPfad#3\endSPfad}
\def\ZeichneSPfad#1{\ifx#1\endSPfad\let\next\relax
  \else\let\next\ZeichneSPfad
    \ifnum#1=1
      \raise\ycoord \Einheit\hbox to0pt{\hskip\xcoord \Einheit
         \hSSchritt\hss}%
      \advance\xcoord by 1
    \else\ifnum#1=2
      \raise\ycoord \Einheit\hbox to0pt{\hskip\xcoord \Einheit
        \hbox{\hskip-2pt \vSSchritt}\hss}%
      \advance\ycoord by 1
    \else\ifnum#1=3
      \raise\ycoord \Einheit\hbox to0pt{\hskip\xcoord \Einheit
         \DSSchritt\hss}
      \advance\xcoord by 1
      \advance\ycoord by 1
    \else\ifnum#1=4
      \raise\ycoord \Einheit\hbox to0pt{\hskip\xcoord \Einheit
         \dSSchritt\hss}
      \advance\xcoord by 1
      \advance\ycoord by -1
    \fi\fi\fi\fi
  \fi\next}
\def\Koordinatenachsen(#1,#2){\unskip
 \hbox to0pt{\hskip-.5pt\vrule height#2 \Einheit width.5pt depth1 \Einheit}%
 \hbox to0pt{\hskip-1 \Einheit \xcoord#1 \advance\xcoord by1
    \vrule height0.25pt width\xcoord \Einheit depth0.25pt\hss}}
\def\Koordinatenachsen(#1,#2)(#3,#4){\unskip
 \hbox to0pt{\hskip-.5pt \ycoord-#4 \advance\ycoord by1
    \vrule height#2 \Einheit width.5pt depth\ycoord \Einheit}%
 \hbox to0pt{\hskip-1 \Einheit \hskip#3\Einheit 
    \xcoord#1 \advance\xcoord by1 \advance\xcoord by-#3 
    \vrule height0.25pt width\xcoord \Einheit depth0.25pt\hss}}
\def\Gitter(#1,#2){\unskip \xcoord0 \ycoord0 \leavevmode
  \LOOP\ifnum\ycoord<#2
    \loop\ifnum\xcoord<#1
      \raise\ycoord \Einheit\hbox to0pt{\hskip\xcoord \Einheit\Punkt\hss}%
      \advance\xcoord by1
    \repeat
    \xcoord0
    \advance\ycoord by1
  \REPEAT}
\def\Gitter(#1,#2)(#3,#4){\unskip \xcoord#3 \ycoord#4 \leavevmode
  \LOOP\ifnum\ycoord<#2
    \loop\ifnum\xcoord<#1
      \raise\ycoord \Einheit\hbox to0pt{\hskip\xcoord \Einheit\Punkt\hss}%
      \advance\xcoord by1
    \repeat
    \xcoord#3
    \advance\ycoord by1
  \REPEAT}
\def\Label#1#2(#3,#4){\unskip \xdim#3 \Einheit \ydim#4 \Einheit
  \def\lo{\advance\xdim by-.5 \Einheit \advance\ydim by.5 \Einheit}%
  \def\llo{\advance\xdim by-.25cm \advance\ydim by.5 \Einheit}%
  \def\loo{\advance\xdim by-.5 \Einheit \advance\ydim by.25cm}%
  \def\o{\advance\ydim by.25cm}%
  \def\ro{\advance\xdim by.5 \Einheit \advance\ydim by.5 \Einheit}%
  \def\rro{\advance\xdim by.25cm \advance\ydim by.5 \Einheit}%
  \def\roo{\advance\xdim by.5 \Einheit \advance\ydim by.25cm}%
  \def\l{\advance\xdim by-.30cm}%
  \def\r{\advance\xdim by.30cm}%
  \def\lu{\advance\xdim by-.5 \Einheit \advance\ydim by-.6 \Einheit}%
  \def\llu{\advance\xdim by-.25cm \advance\ydim by-.6 \Einheit}%
  \def\luu{\advance\xdim by-.5 \Einheit \advance\ydim by-.30cm}%
  \def\u{\advance\ydim by-.30cm}%
  \def\ru{\advance\xdim by.5 \Einheit \advance\ydim by-.6 \Einheit}%
  \def\rru{\advance\xdim by.25cm \advance\ydim by-.6 \Einheit}%
  \def\ruu{\advance\xdim by.5 \Einheit \advance\ydim by-.30cm}%
  #1\raise\ydim\hbox to0pt{\hskip\xdim
     \vbox to0pt{\vss\hbox to0pt{\hss$#2$\hss}\vss}\hss}%
}
\catcode`\@=13

\magnification1200
\hsize13cm
\vsize19cm

\TagsOnRight

\def\StemAE{11}
\def\SagaAL{10}
\def\PropAA{9}
\def\MiRRAB{8}
\def\MacdAC{7}
\def\HeGeAA{6}
\def\GeViAB{5}
\def\FuKrAA{4}
\def\FulmAB{3}
\def\ElKLAA{2}
\def\CiucAB{1}

\def\fl#1{\left\lfloor#1\right\rfloor}
\def\cl#1{\left\lceil#1\right\rceil}
\def\v#1{\left\vert#1\right\vert}
\def\prodl{\prod\limits}
\def\la{\lambda}
\def\al{\alpha}
\def\be{\beta}

\topmatter 
\title Schur function identities and the number of perfect matchings of
holey Aztec rectangles
\endtitle 
\author C.~Krattenthaler
\endauthor 
\affil 
Institut f\"ur Mathematik der Universit\"at Wien,\\
Strudlhofgasse 4, A-1090 Wien, Austria.\\
e-mail: KRATT\@Pap.Univie.Ac.At\\
WWW: \tt http://radon.mat.univie.ac.at/People/kratt
\endaffil
\address Institut f\"ur Mathematik der Universit\"at Wien,
Strudlhofgasse 4, A-1090 Wien, Austria.
\endaddress
\abstract We compute the number of perfect matchings of an $M\times N$ Aztec
rectangle where $\vert N-M\vert$ vertices have been removed along a line. A
particular case solves a problem posed by Propp. 
Our enumeration results follow from
certain identities for Schur functions, which
are established by the combinatorics of nonintersecting lattice paths.
\endabstract
\endtopmatter
\document

\rightheadtext{Schur function identities and the number of perfect matchings}

\subhead 1. Introduction\endsubhead
Consider a $(2M+1)\times (2N+1)$ rectangular chessboard and suppose
that the corners are black. The $M\times N$ {\it Aztec rectangle} is
the graph whose vertices are the white squares and whose edges
connect precisely those pairs of white squares that are diagonally
adjacent, see Figure~1.a for an example. An $M\times N$ Aztec rectangle
does not have any perfect matching, except if $M=N$. If, however, we
remove $\vert N-M\vert$ vertices from it in a suitable way, then the
``holey" Aztec rectangle which is obtained in this manner 
(see Figures~4 and 5 for examples)
does allow perfect matchings. For example, Ciucu \cite{\CiucAB,
Theorem~4.1} computed, for even $M$ and $M\le N$, the precise
number of all perfect matchings of an
$M\times N$ Aztec rectangle where $N-M$ vertices on the horizontal
symmetry axis are removed (see Theorem~7).

In Problem~10 of his list \cite{\PropAA} of ``20 open problems on
enumeration of matchings", Propp asks for the number 
of all perfect matchings of a
$(2n-1)\times 2n$ Aztec rectangle where one vertex which is
adjacent to the central vertex is removed (it does not matter which
of the four vertices adjacent to the central one is removed; see
Figure~1.b for an example).
\vskip10pt
\vbox{
$$
\Einheit.45cm
\PfadDicke{.5pt}
\thinlines
\Pfad(0,1),33333444\endPfad
\Pfad(0,3),33344444\endPfad
\Pfad(0,5),34444443\endPfad
\Pfad(0,1),43333334\endPfad
\Pfad(0,3),44433333\endPfad
\Pfad(0,5),44444333\endPfad
\DickPunkt(0,1)
\DickPunkt(0,3)
\DickPunkt(0,5)
\DickPunkt(1,0)
\DickPunkt(1,2)
\DickPunkt(1,4)
\DickPunkt(1,6)
\DickPunkt(2,1)
\DickPunkt(2,3)
\DickPunkt(2,5)
\DickPunkt(3,0)
\DickPunkt(3,2)
\DickPunkt(3,4)
\DickPunkt(3,6)
\DickPunkt(4,1)
\DickPunkt(4,3)
\DickPunkt(4,5)
\DickPunkt(5,0)
\DickPunkt(5,2)
\DickPunkt(5,4)
\DickPunkt(5,6)
\DickPunkt(6,1)
\DickPunkt(6,3)
\DickPunkt(6,5)
\DickPunkt(7,0)
\DickPunkt(7,2)
\DickPunkt(7,4)
\DickPunkt(7,6)
\DickPunkt(8,1)
\DickPunkt(8,3)
\DickPunkt(8,5)
\hbox{\hskip5cm}
\Pfad(0,1),33333444\endPfad
\Pfad(0,3),33344444\endPfad
\Pfad(7,2),4\endPfad
\Pfad(0,5),3444\endPfad
\Pfad(6,1),43\endPfad
\Pfad(0,1),43333334\endPfad
\Pfad(0,3),4443\endPfad
\Pfad(6,3),33\endPfad
\Pfad(0,5),44444333\endPfad
\DickPunkt(0,1)
\DickPunkt(0,3)
\DickPunkt(0,5)
\DickPunkt(1,0)
\DickPunkt(1,2)
\DickPunkt(1,4)
\DickPunkt(1,6)
\DickPunkt(2,1)
\DickPunkt(2,3)
\DickPunkt(2,5)
\DickPunkt(3,0)
\DickPunkt(3,2)
\DickPunkt(3,4)
\DickPunkt(3,6)
\DickPunkt(4,1)
\DickPunkt(4,3)
\DickPunkt(4,5)
\DickPunkt(5,0)
\Kreis(5,2)
\DickPunkt(5,4)
\DickPunkt(5,6)
\DickPunkt(6,1)
\DickPunkt(6,3)
\DickPunkt(6,5)
\DickPunkt(7,0)
\DickPunkt(7,2)
\DickPunkt(7,4)
\DickPunkt(7,6)
\DickPunkt(8,1)
\DickPunkt(8,3)
\DickPunkt(8,5)
\hskip4.3cm
$$
\centerline{\eightpoint a. A $3\times 4$ Aztec rectangle\hskip1cm
b. A $3\times 4$ Aztec rectangle with a vertex}
\centerline{\hskip5cm\eightpoint 
adjacent to the central vertex removed}
\vskip10pt
\centerline{\eightpoint Figure 1}
}
\vskip10pt

The original purpose of this paper was to solve this problem. (An independent
solution was found by Helfgott and Gessel \cite{\HeGeAA, Sec.~6}.)
What we do in fact is, 
quite generally, to describe a method which allows to compute the 
number of all perfect matchings of an
$M\times N$ Aztec rectangle where $\vert N-M\vert$ vertices are
removed along a horizontal line. The expressions that are obtained
(see Theorems~11 and 12)
are $\fl{d/2}$-fold sums, where $d$ is the ``distance" of the line
from the symmetry axis (``distance"
has the obvious meaning here, see Section~4). Thus, the closer the
horizontal line is to the symmetry axis, the simpler are the obtained
expressions. In particular, closed forms are obtained if the line is
the symmetry axis or ``by 1 off" the symmetry axis. This includes 
Ciucu's enumeration and the case of Propp's problem (see
Theorems~7--10). In addition, with the help of hypergeometric-type summations
(Theorem~6), we also obtain closed forms if the vertices
which are not removed ``form an arithmetic progression" (Theorems~13 and 14) 
or ``form a geometric progression" (Theorems~15 and 16), 
regardless on which line they are removed.

As a basis for our method we take a careful selection of ideas which
appear in the literature in connection with this theme. To be more
concrete, Ciucu proves his above cited enumeration result by applying his
matchings factorization theorem for symmetric graphs \cite{\CiucAB,
Theorem~1.2} and then using a
formula due to Mills, Robbins and Rumsey \cite{\MiRRAB, Theorem~2}
(see Lemma~1)
for the number of perfect matchings of an $M\times N$ Aztec rectangle
with $\vert N-M\vert$ of its top-most vertices removed. This approach
does of course not generalize since there is no symmetry if vertices
are removed arbitrarily along a horizontal line {\it different\/}
from the symmetry axis. However, as an afterthought \cite{\CiucAB,
Remark~4.3}, he sketches a second approach, which again makes use of
the formula due to Mills, Robbins and Rumsey, and which finally boils
down to establishing a certain Schur function identity (see Theorem~3),
which he states as a conjecture. It is this approach that we are
going to adapt (in a slightly modified fashion). 
The conjectured Schur function identity was subsequently
proved by Tesler (private communication) making use of Laplace expansion
of a certain matrix, and independently by Fulmek \cite{\FulmAB}
making use of nonintersecting lattice paths. Whereas Tesler's proof
does not seem to be of any help for the generalizations, Fulmek's
idea of proof is exactly the right tool for proving the variations of
this Schur function identity that we need.

Our paper is organized as follows. In the next section we list two
auxiliary results, one of which is the formula of Mills, Robbins and
Rumsey and the other is a variation of it. In Section~3 we review Fulmek's lattice
path proof of Ciucu's (conjectured) Schur function identity and
describe how the same idea leads to a whole family of variations of
the identity. Also contained in this section is the
hypergeometric-type summation that was mentioned earlier.
Finally, in Section~4, we describe our method how to enumerate 
all perfect matchings of an
$M\times N$ Aztec rectangle where $\vert N-M\vert$ vertices are
removed along a horizontal line, and list our results.
The method simply consists of breaking
the Aztec rectangle into two parts along this horizontal line, and
then handling the resulting summations by taking advantage of the
Schur function identities from Section~3. 
The impatient reader may immediately jump to this section, and come
back to Section~2, respectively Section~3, when results from there are cited.

\subhead 2. Two auxiliary lemmas\endsubhead
First, we quote Theorem~2 from \cite{\MiRRAB}, in the translation as
described by
Ciucu \cite{\CiucAB, (4.4)} (see also equation (7) in \cite{\ElKLAA}).
\proclaim{Lemma~1}
The number of perfect matchings of an $m\times n$
Aztec rectangle, where all the vertices in the top-most row, except for 
the $a_1$-st, the $a_2$-nd, \dots, and the $a_{m}$-th vertex, have been
removed (see Figure~2.a for an example with $m=3$, $n=5$, $a_1=1$,
$a_2=3$, $a_3=5$), equals
$$\frac {2^{\binom {m+1}2}} {\prodl _{i=1} ^{m}(i-1)!}\prodl _{1\le
i<j\le m} ^{}(a_j-a_i).\tag2.1$$
\endproclaim
\midinsert
\vskip10pt
\vbox{
$$
\Einheit.45cm
\PfadDicke{.5pt}
\thinlines
\Pfad(0,1),3333344444\endPfad
\Pfad(0,3),33\endPfad
\Pfad(4,5),444443\endPfad
\Pfad(3,4),444433\endPfad
\Pfad(0,5),3444444333\endPfad
\Pfad(0,1),433333\endPfad
\Pfad(8,5),44\endPfad
\Pfad(0,3),4443333334\endPfad
\Pfad(0,5),4444433333\endPfad
\DickPunkt(0,1)
\DickPunkt(0,3)
\DickPunkt(0,5)
\DickPunkt(1,0)
\DickPunkt(1,2)
\DickPunkt(1,4)
\DickPunkt(1,6)
\DickPunkt(2,1)
\DickPunkt(2,3)
\DickPunkt(2,5)
\DickPunkt(3,0)
\DickPunkt(3,2)
\DickPunkt(3,4)
\Kreis(3,6)
\DickPunkt(4,1)
\DickPunkt(4,3)
\DickPunkt(4,5)
\DickPunkt(5,0)
\DickPunkt(5,2)
\DickPunkt(5,4)
\DickPunkt(5,6)
\DickPunkt(6,1)
\DickPunkt(6,3)
\DickPunkt(6,5)
\DickPunkt(7,0)
\DickPunkt(7,2)
\DickPunkt(7,4)
\Kreis(7,6)
\DickPunkt(8,1)
\DickPunkt(8,3)
\DickPunkt(8,5)
\DickPunkt(9,0)
\DickPunkt(9,2)
\DickPunkt(9,4)
\DickPunkt(9,6)
\DickPunkt(10,1)
\DickPunkt(10,3)
\DickPunkt(10,5)
\hbox{\hskip6.5cm}
\Pfad(1,0),333334444\endPfad
\Pfad(0,1),333\endPfad
\Pfad(5,4),44443\endPfad
\Pfad(0,3),3344444333\endPfad
\Pfad(0,5),44444333\endPfad
\Pfad(3,0),3333\endPfad
\Pfad(9,4),4\endPfad
\Pfad(0,1),4\endPfad
\Pfad(3,0),3333\endPfad
\Pfad(0,3),444\endPfad
\Pfad(5,0),3333\endPfad
\DickPunkt(0,1)
\DickPunkt(0,3)
\DickPunkt(0,5)
\DickPunkt(1,0)
\DickPunkt(1,2)
\DickPunkt(1,4)
\DickPunkt(2,1)
\DickPunkt(2,3)
\DickPunkt(2,5)
\DickPunkt(3,0)
\DickPunkt(3,2)
\DickPunkt(3,4)
\DickPunkt(4,1)
\DickPunkt(4,3)
\Kreis(4,5)
\DickPunkt(5,0)
\DickPunkt(5,2)
\DickPunkt(5,4)
\DickPunkt(6,1)
\DickPunkt(6,3)
\DickPunkt(6,5)
\DickPunkt(7,0)
\DickPunkt(7,2)
\DickPunkt(7,4)
\DickPunkt(8,1)
\DickPunkt(8,3)
\Kreis(8,5)
\DickPunkt(9,0)
\DickPunkt(9,2)
\DickPunkt(9,4)
\DickPunkt(10,1)
\DickPunkt(10,3)
\Kreis(10,5)
\Kreis(1,6)
\Kreis(3,6)
\Kreis(5,6)
\Kreis(7,6)
\Kreis(9,6)
\hskip4.5cm
$$
\line{\eightpoint\hskip.5cm a. A $3\times 5$ Aztec rectangle with some vertices\hskip.5cm
b. A $3\times 5$ Aztec rectangle with all vertices\hss}
\line{\eightpoint \hphantom{\hskip.5cm a. }in the top row removed \hskip3.35cm 
in the top row removed and some\hss}
\line{\eightpoint 
\hphantom{\hskip.5cm a. A $3\times 5$ Aztec rectangle with some vertices\hskip.5cm
b. }vertices in the next row removed\hss}
\vskip10pt
\centerline{\eightpoint Figure 2}
}
\vskip10pt
\endinsert
Next, we state a variant of this result for an Aztec rectangle with
a holey ``half-row" (see the last paragraph of Part~1 of \cite{\ElKLAA},
and also \cite{\HeGeAA, Lemma~2}).
\proclaim{Lemma~2}
The number of perfect matchings of an $m\times n$
Aztec rectangle, where all the vertices in the top-most row have been
removed, and where 
the $a_1$-st, the $a_2$-nd, \dots, and the $a_{m}$-th vertex of the
row next to the top-most row have been
removed (see Figure~2.b for an example with $m=3$, $n=5$, $a_1=3$,
$a_2=5$, $a_3=6$), equals
$$\frac {2^{\binom {m}2}} {\prodl _{i=1} ^{m}(i-1)!}\prodl _{1\le
i<j\le m} ^{}(a_j-a_i).\tag2.2$$
\endproclaim

\subhead 3. Schur function identities\endsubhead
Let $A=\{a_1,a_2,\dots,a_m\}$ be a set of positive integers with
$a_1< a_2<\dots<a_m$.
The connection between enumeration of perfect matchings of holey
Aztec rectangles and Schur functions is set up by the identity (see
\cite{\MacdAC, Ex.~I.3.4})
$$s_{\la(A)}(1^m)=s_{\la(A)}(\underbrace{1,1,\dots,1}_{m\text {
times}})=\frac {\prodl _{1\le i<j\le m} ^{}(a_j-a_i)} {\prodl _{i=1}
^{m}(i-1)!},\tag3.1$$
where $\la(A)$ is by definition the partition
$(a_m-m+1,\dots,a_2-1,a_1)$. It is the ubiquitous ``Vandermonde"
product, as it appears in Lemmas~1 and 2 and (3.1), which makes it
natural to consider Schur functions in this context, see Section~4.

Let us recall Ciucu's (conjectured) Schur function 
identity \cite{\CiucAB, Conj.~4.5} and Fulmek's \cite{\FulmAB} proof
of it.
Here and in the sequel, we use the short notation $X_n$ for the
sequence of variables $x_1,x_2,\dots,x_n$.
\proclaim{Theorem~3}Let $T=\{t_1,t_2,\dots,t_{2m}\}$ be a set of
positive integers with $t_1<t_2<\dots<t_{2m}$. Then
$$\sum _{} ^{}s_{\la(A)}(X_n)\cdot s_{\la(B)}(X_n)=2^m
s_{\la(t_2,t_4,\dots,t_{2m})}(X_n)\cdot 
s_{\la(t_1,t_3,\dots,t_{2m-1})}(X_n),\tag3.2$$
where the sum is over all pairs of disjoint sets $A$ and $B$ whose
union is $T$ and whose cardinalities are given by $\v{A}=\v{B}=m$.
\endproclaim
\demo{Sketch of Proof} By the main theorem of nonintersecting lattice paths 
\cite{\GeViAB, Cor.~2; \StemAE, Theorem~1.2}, a Schur function
$s_{\la(A)}(X_n)$ can be interpreted combinatorially as the
generating function $\sum _{\Cal P} ^{}w(\Cal P)$, where the sum is
over all families $\Cal P=(P_1,P_2,\dots,P_m)$ of nonintersecting
lattice paths consisting of horizontal and vertical unit steps, 
$P_i$ running from $(i-1,1)$ to $(a_i,n)$,
$i=1,2,\dots,m$, and where the weight $w(.)$ is defined via edge
weights in which vertical edges have weight 1 and a horizontal edge
at height $h$ has weight $x_h$ (see also \cite{\FuKrAA, Sec.~3;
\SagaAL, Sec.~4.5}).

Thus, the left-hand side of (3.2) can be interpreted as a generating
function for pairs $(\Cal P^g,\Cal P^r)$ of families of nonintersecting 
lattice paths, the $i$-th path of $\Cal P^g$ 
running from $(i-1,1)$ to $(a_i,n)$, with $a_i$ being the $i$-th
element of $A$, $i=1,2,\dots,m$, and the $i$-th path of $\Cal P^r$ 
running from $(i-1,1)$ to $(b_i,n)$, with $b_i$ being the $i$-th
element of $B$, $i=1,2,\dots,m$. 
Say that the paths of $\Cal P^g$ are coloured green,
and those of $\Cal P^r$ are coloured red. Likewise, the 
right-hand side of (3.2) can be interpreted as a generating
function for triples $(\Cal Q^g,\Cal Q^r,C)$, where $\Cal Q^g$ and
$\Cal Q^r$ are families of nonintersecting 
lattice paths, the $i$-th path of $\Cal Q^g$ 
running from $(i-1,1)$ to $(t_{2i},n)$, $i=1,2,\dots,m$, 
and the $i$-th path of $\Cal Q^r$ 
running from $(i-1,1)$ to $(t_{2i-1},n)$, $i=1,2,\dots,m$, 
and where $C$ is a $\{0,1\}$-sequence
of length $m$.

Identity (3.2) will be proved once we are able to set up a
(weight-preserving) bijection between ``left-hand side pairs" and
``right-hand side triples". Fulmek does this in the following way.
Let $(\Cal P^g,\Cal P^r)$ be a ``left-hand side pair" as described
above. Then Fulmek defines a (non-crossing) matching of the end
points $(t_i,n)$, $i=1,2,\dots,2m$, by ``down-up" trails. Beginning
in some end point we move {\it down} along the path which ends in
this point. When we meet a path of colour different from the path
along which we were moving, then we continue
to move {\it up} along the new path. This procedure is iterated,
interchanging the roles of up and down every time, of
course. Finally, we will terminate, in an up-move, in another end
point. It is easy to see that in that manner an end point with odd
index, $(t_{2i-1},n)$ say, is always connected with an end point of
even index, $(t_{2j},n)$ say. The triple $(\Cal Q^g,\Cal Q^r,C)$ is
now defined as follows. For any $i$, $i=1,2,\dots,m$, we consider the
down-up trail starting in $(t_{2i-1},n)$. If the path ending in that
end point is green then we interchange colours along the down-up trail
and we put $C_i=1$. On the other hand, if the path ending in that
end point is red then we leave the paths as they are
and put $C_i=0$. By definition, $\Cal Q^g$ is the family of green paths thus
obtained, $\Cal Q^r$ is the family of red paths thus
obtained, and $C$ is the sequence $C_1,C_2,\dots,C_m$ thus obtained.

It is easy to see that this mapping has all the required properties.
We refer the reader to \cite{\FulmAB}.\quad \quad \qed
\enddemo
For the solution of Propp's Problem~10, we need a variant of
Theorem~3.
\proclaim{Theorem~4}Let $T=\{t_1,t_2,\dots,t_{2m+1}\}$ be a set of
positive integers with $t_1<t_2<\dots<t_{2m+1}$. Then
$$\sum _{} ^{}s_{\la(A)}(X_n)\cdot s_{\la(B)}(X_{n+1})=2^m
s_{\la(t_2,t_4,\dots,t_{2m})}(X_n)\cdot 
s_{\la(t_1,t_3,\dots,t_{2m+1})}(X_{n+1}),\tag3.3$$
where the sum is over all disjoint pairs of sets $A$ and $B$ whose
union is $T$ and whose cardinalities are given by $\v{A}=m$ and
$\v{B}=m+1$.
\endproclaim
\demo{Sketch of Proof} We proceed in the same manner as in the proof
of Theorem~3. In particular, the left-hand side of (3.3) can be 
interpreted as a generating
function for pairs $(\Cal P^g,\Cal P^r)$ of families of nonintersecting 
lattice paths, the $i$-th path of $\Cal P^g$ 
running from $(i-1,1)$ to $(a_i,n)$, with $a_i$ being the $i$-th
element of $A$, $i=1,2,\dots,m$, and the $i$-th path of $\Cal P^r$ 
running from $(i-1,0)$ to $(b_i,n)$, with $b_i$ being the $i$-th
element of $B$, $i=1,2,\dots,m+1$. 
Likewise, the 
right-hand side of (3.3) can be interpreted as a generating
function for triples $(\Cal Q^g,\Cal Q^r,C)$, where $\Cal Q^g$ and
$\Cal Q^r$ are families of nonintersecting 
lattice paths, the $i$-th path of $\Cal Q^g$ 
running from $(i-1,1)$ to $(t_{2i},n)$, $i=1,2,\dots,m$, 
and the $i$-th path of $\Cal Q^r$ 
running from $(i-1,0)$ to $(t_{2i-1},n)$, $i=1,2,\dots,m+1$, 
and where $C$ is a $\{0,1\}$-sequence
of length $m$.

For setting up a bijection between ``left-hand side pairs" and
``right-hand side triples", we construct again the down-up trails for
a pair $(\Cal P^g,\Cal P^r)$.
However, since the number of end points is now odd, it is impossible
to obtain a complete matching of the end points. In fact, what
is obtained is a (non-crossing) matching on only $2m$ end points,
while a single end point is ``matched" to the
``additional" starting point of red paths, $(m,0)$. Besides, the
corresponding down-up trail, which connects this single end point with
$(m,0)$, cannot touch any of the other down-up trails. Hence, it has
odd index, i.e., it is an end point $(t_{2\ell-1},n)$, say.
Furthermore, since this trail must end with a ``down part", moving
into $(m,0)$, which is the starting point of a red path only, and since
all ``down parts" are equally coloured red (as well as all the ``up parts"
are coloured in the opposite colour green), the end point
$(t_{2\ell-1},n)$ is the end point of a red path.

To define the triple $(\Cal Q^g,\Cal Q^r,C)$ we must modify the
definition in the proof of Theorem~3 only slightly. 
Namely, instead of considering the
down-up trail starting in $(t_{2i-1},n)$ for {\it any} $i$ between
$1$ and $m+1$, and possibly performing a recolouring and 
defining $C_i$, we do
it only for all $i$ between $1$ and $m+1$ {\it that are different from
$\ell$}. (Recall that $(t_{2\ell-1},n)$ is matched to a starting
point and that the path ending in this point must be red.) 
Everything else remains the same, except that the sequence $C$ is now
defined as $C_1,C_2,\dots,C_{\ell-1},C_{\ell+1},\dots,C_{m+1}$.\quad \quad
\qed
\enddemo

Clearly, there is nothing which could prevent us from considering the
next summation in this family, of which (3.2) and (3.3) are the first
members,
$$\sum _{} ^{}s_{\la(A)}(X_n)\cdot s_{\la(B)}(X_{n+2}),$$
where the sum is over all disjoint pairs of sets $A$ and $B$ whose
union is $T=\{t_1,t_2,\dots,\mathbreak t_{2m+2}\}$ 
and whose cardinalities are given by $\v{A}=m$ and
$\v{B}=m+2$. In general, we would consider
$$\sum _{} ^{}s_{\la(A)}(X_n)\cdot s_{\la(B)}(X_{n+d}),$$
where the sum is over all disjoint pairs of sets $A$ and $B$ whose
union is $T=\{t_1,t_2,\dots,\mathbreak t_{2m+d}\}$ 
and whose cardinalities are given by $\v{A}=m$ and
$\v{B}=m+d$. Of course, the matching and ``colouring scheme" will not
lead to closed forms in general, the reason being that $d$ end points
will be matched to starting points. What is obtained is a result
in form of a $\fl{d/2}$-fold sum, which we state next. 

\proclaim{Theorem~5}Let $T=\{t_1,t_2,\dots,t_{2m+d}\}$ be a set of
positive integers with $t_1<t_2<\dots<t_{2m+d}$. Then
$$\multline
\sum _{} ^{}s_{\la(A)}(X_n)\cdot s_{\la(B)}(X_{n+d})\\
=2^m
\sum _{1\le k_1<k_2<\dots<k_{\fl{d/2}}\le m+\fl{d/2}} ^{}
s_{\la(\{t_2,t_4,\dots,t_{2m+d}\}\backslash
\{2k_1,2k_2,\dots,2k_{\fl{d/2}}\})}
(X_n)\\
\cdot 
s_{\la(\{t_1,t_3,\dots,t_{2m+1}\}\cup 
\{2k_1,2k_2,\dots,2k_{\fl{d/2}}\})}(X_{n+d}),
\endmultline\tag3.4$$
where the sum is over all disjoint pairs of sets $A$ and $B$ whose
union is $T$ and whose cardinalities are given by $\v{A}=m$ and
$\v{B}=m+d$, and where complement and union in the indices on the
right-hand side have the obvious meaning. (In particular, in case of
`union' the elements are brought in order before applying the
$\la$-operation).
\endproclaim

The last piece that we need  for our enumerations
is summations for sums that will result
from application of Theorem~5. Although they are (basic)
hypergeometric-type multiple sums, they are again actually a consequence of a
Schur function identity.
\proclaim{Theorem~6}With the usual definition of shifted factorials, 
$(a)_k:=a(a+1)\cdots(a+k-1)$ for
$k\ge1$, and $(a)_0:=1$, we have
$$\multline \sum _{0\le k_1<k_2<\dots<k_s\le m} ^{}\prodl _{1\le
i<j\le s} ^{}(k_j-k_i)^2\prodl _{i=1} ^{s}\frac
{(x)_{k_i}\,(y)_{m-k_i}} {(k_i)!\,(m-k_i)!}\\
=\prodl _{i=1} ^{s}\frac
{(x)_{i-1}\,(y)_{i-1}\,(x+y+i+s-2)_{m-s+1}\,(i-1)!} {(m-i+1)!}.
\endmultline\tag3.5$$ 
With the usual definition of shifted $q$-factorials, 
$(a;q)_k:=(1-a)(1-aq)\cdots(1-aq^{k-1})$ for
$k\ge1$, and $(a;q)_0:=1$, we have
$$\multline \sum _{0\le k_1<k_2<\dots<k_s\le m} ^{}y^{\sum _{i=1}
^{s}k_i}\prodl _{1\le
i<j\le s} ^{}(q^{k_j}-q^{k_i})^2\prodl _{i=1} ^{s}\frac
{(x;q)_{k_i}\,(y;q)_{m-k_i}} {(q;q)_{k_i}\,(q;q)_{m-k_i}}\\
=q^{\binom s3}y^{\binom s2}\prodl _{i=1} ^{s}\frac
{(x;q)_{i-1}\,(y;q)_{i-1}\,(xyq^{i+s-2};q)_{m-s+1}\,(q;q)_{i-1}}
{(q;q)_{m-i+1}}.
\endmultline\tag3.6$$ 
\endproclaim
\demo{Proof}It suffices to prove identity (3.6). For, identity (3.5)
immediately results from (3.6), by performing the replacements $x\to
q^x$ and $y\to q^y$ in (3.6), then dividing both sides by
$(1-q)^{s(s-1)}$, and subsequently performing the limit $q\to 1$.

We start with the well-known Schur function identity
(see \cite{\MacdAC, I, (5.9)})
$$s_\la(x_1,\dots,x_\al ,y_1,\dots,y_\be)=\sum _{\mu}
^{}s_{\la/\mu}(x_1,\dots,x_\al ) s_{\mu}(y_1,\dots,y_\be),\tag3.7$$
where the sum is over all partitions $\mu$ that are contained in the
partition $\la$. The special case which is relevant for us is the
case that $\la$ is a rectangle, $\la=(M^s)$ say. Then the skew Schur
function $s_{\la/\mu}(x_1,\dots,x_\al )=s_{(M^s)/\mu}(x_1,\dots,x_\al )$ 
equals an ordinary Schur
function, namely we have
$$s_{(M^s)/\mu}(x_1,\dots,x_\al )=s_{\mu'}(x_1,\dots,x_\al ),\tag3.8$$
where, given $\mu=(\mu_1,\mu_2,\dots,\mu_s)$, the partition $\mu'$ is
defined by $\mu'=(M-\mu_s,M-\mu_{s-1},\dots,M-\mu_1)$.

Now we set $x_i=q^i$, $i=1,2,\dots,\al$, and $y_i=q^{\al+i}$,
$i=1,2,\dots,\be$, in (3.7). Use of
(3.8),
of the formula (see \cite{\MacdAC, ??})
$$\multline
s_\mu(q^{K+1},q^{K+2},\dots,q^{K+L})\\
=q^{2\binom {s+1}3+(K+1)\sum _{i=1} ^{s}\mu_i}\frac {\prodl _{1\le i<j\le s} ^{}(q^{\mu_j-j}-q^{\mu_i-i}) 
\prodl _{i=1} ^{s}(q^{L+1-i};q)_{\mu_i}} {\prodl _{i=1}
^{s}(q;q)_{\mu_i-i+s}},
\endmultline$$
and routine manipulations show, 
that (3.7) is equivalent to (3.6), with the relations
$k_i=\mu_i-i+s$, $x=q^{\al-s+1}$, $y=q^{\be-s+1}$, $m=M+s-1$.

Hence, we know that (3.6) is true if $x$ and $y$ are positive integer powers
of $q$.
In order to extend (3.6) to all $x$ and $y$, we observe that, first, for
fixed $s$ and $m$, left-hand side and right-hand side of (3.6) are
polynomials in $x$ and $y$ with degree bounded by a fixed quantity, 
and second, that (3.6) is true for an infinite number of $x$'s and
$y$'s. Therefore (3.6) must be true in general.\quad \quad \qed
\enddemo

\subhead 4. The enumeration results\endsubhead
Now we are ready to state our enumeration results for holey Aztec
rectangles. For the sake of completeness, we start by repeating
Ciucu's enumeration that was mentioned in the Introduction. 
\proclaim{Theorem~7 \rm{\cite{\CiucAB, Theorem~4.1}}}
Let $m$ and $N$ be positive integers with $2m\le N$. Then 
the number of perfect matchings of a $2m\times N$
Aztec rectangle, where all the vertices on the central horizontal
row, except for 
the $t_1$-st, the $t_2$-nd, \dots, and the $t_{2m}$-th vertex, have been
removed (see Figure~4.a for an example with $m=2$, $N=7$, $a_1=1$,
$a_2=4$, $a_3=5$, $a_4=7$), equals
$$\frac {2^{m^2+2m}} {\prodl _{i=1} ^{m}(i-1)!^2}\prodl _{1\le
i<j\le m} ^{}(t_{2j}-t_{2i})\prodl _{1\le
i<j\le m} ^{}(t_{2j-1}-t_{2i-1}).\tag4.1$$
\endproclaim
\vskip10pt
\vbox{
$$
\Einheit.4cm
\PfadDicke{.5pt}
\thinlines
\Pfad(0,3),43434343434343\endPfad
\Pfad(0,3),34343434343434\endPfad
\Pfad(0,1),43\endPfad
\Pfad(6,1),4343\endPfad
\Pfad(12,1),43\endPfad
\Pfad(0,1),34343434343434\endPfad
\Pfad(0,-1),43434343434343\endPfad
\Pfad(0,-1),34\endPfad
\Pfad(6,-1),3434\endPfad
\Pfad(12,-1),34\endPfad
\Pfad(0,-3),43434343434343\endPfad
\Pfad(0,-3),34343434343434\endPfad
\DickPunkt(1,0)
\Kreis(3,0)
\Kreis(5,0)
\DickPunkt(7,0)
\DickPunkt(9,0)
\Kreis(11,0)
\DickPunkt(13,0)
\DickPunkt(0,1)
\DickPunkt(2,1)
\DickPunkt(4,1)
\DickPunkt(6,1)
\DickPunkt(8,1)
\DickPunkt(10,1)
\DickPunkt(12,1)
\DickPunkt(14,1)
\DickPunkt(1,2)
\DickPunkt(3,2)
\DickPunkt(5,2)
\DickPunkt(7,2)
\DickPunkt(9,2)
\DickPunkt(11,2)
\DickPunkt(13,2)
\DickPunkt(0,3)
\DickPunkt(2,3)
\DickPunkt(4,3)
\DickPunkt(6,3)
\DickPunkt(8,3)
\DickPunkt(10,3)
\DickPunkt(12,3)
\DickPunkt(14,3)
\DickPunkt(1,4)
\DickPunkt(3,4)
\DickPunkt(5,4)
\DickPunkt(7,4)
\DickPunkt(9,4)
\DickPunkt(11,4)
\DickPunkt(13,4)
\DickPunkt(1,-2)
\DickPunkt(3,-2)
\DickPunkt(5,-2)
\DickPunkt(7,-2)
\DickPunkt(9,-2)
\DickPunkt(11,-2)
\DickPunkt(13,-2)
\DickPunkt(0,-1)
\DickPunkt(2,-1)
\DickPunkt(4,-1)
\DickPunkt(6,-1)
\DickPunkt(8,-1)
\DickPunkt(10,-1)
\DickPunkt(12,-1)
\DickPunkt(14,-1)
\DickPunkt(1,-4)
\DickPunkt(3,-4)
\DickPunkt(5,-4)
\DickPunkt(7,-4)
\DickPunkt(9,-4)
\DickPunkt(11,-4)
\DickPunkt(13,-4)
\DickPunkt(0,-3)
\DickPunkt(2,-3)
\DickPunkt(4,-3)
\DickPunkt(6,-3)
\DickPunkt(8,-3)
\DickPunkt(10,-3)
\DickPunkt(12,-3)
\DickPunkt(14,-3)
\hbox{\hskip8.5cm}
\Pfad(0,4),434343\endPfad
\Pfad(0,4),343434\endPfad
\Pfad(0,2),434343\endPfad
\Pfad(0,2),343434\endPfad
\Pfad(1,1),43\endPfad
\Pfad(1,-1),34\endPfad
\Pfad(5,1),4\endPfad
\Pfad(5,-1),3\endPfad
\Pfad(0,-2),434343\endPfad
\Pfad(0,-2),343434\endPfad
\Pfad(0,-4),434343\endPfad
\Pfad(0,-4),343434\endPfad
\DickPunkt(1,-1)
\DickPunkt(3,-1)
\DickPunkt(5,-1)
\Kreis(0,0)
\DickPunkt(2,0)
\Kreis(4,0)
\DickPunkt(6,0)
\DickPunkt(1,1)
\DickPunkt(3,1)
\DickPunkt(5,1)
\DickPunkt(0,2)
\DickPunkt(2,2)
\DickPunkt(4,2)
\DickPunkt(6,2)
\DickPunkt(1,3)
\DickPunkt(3,3)
\DickPunkt(5,3)
\DickPunkt(0,4)
\DickPunkt(2,4)
\DickPunkt(4,4)
\DickPunkt(6,4)
\DickPunkt(1,5)
\DickPunkt(3,5)
\DickPunkt(5,5)
\DickPunkt(0,-4)
\DickPunkt(2,-4)
\DickPunkt(4,-4)
\DickPunkt(6,-4)
\DickPunkt(1,-5)
\DickPunkt(3,-5)
\DickPunkt(5,-5)
\DickPunkt(1,-3)
\DickPunkt(3,-3)
\DickPunkt(5,-3)
\DickPunkt(0,-2)
\DickPunkt(2,-2)
\DickPunkt(4,-2)
\DickPunkt(6,-2)
\hskip3.5cm
$$
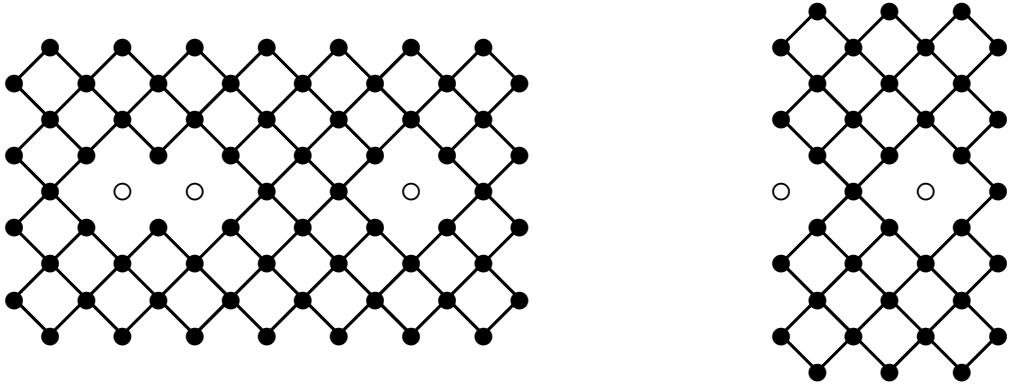
\line{\eightpoint\hskip.1cm a. A $4\times 7$ Aztec rectangle with some 
vertices\hskip.5cm
b. A $5\times 3$ Aztec rectangle with some vertices\hss}
\line{\eightpoint\hphantom{\hskip.1cm a. }on the central row removed\hskip2.95cm
on the central line removed\hss}
\vskip10pt
\centerline{\eightpoint Figure 4}
}
\vskip10pt
Our first new result is an analogue of this theorem for an odd number
of rows. It was also obtained independently by Helfgott and Gessel
\cite{\HeGeAA, Proposition~6}.
\proclaim{Theorem~8}Let $m$ and $N$ be positive integers with $2m-1\ge N$. Then 
the number of perfect matchings of a $(2m-1)\times N$
Aztec rectangle, where all the vertices on the central horizontal
row, except for 
the $t_1$-st, the $t_2$-nd, \dots, and the $t_{2N-2m+2}$-nd vertex,
have been removed (see Figure~4.b for an example with $m=3$, $N=3$,
$a_1=2$, $a_2=4$), equals
$$\multline
{2^{m^2-2m+N+1}}\frac {\prodl _{i=m+1} ^{N+1}(i-1)!^2} {\prodl _{i=1}
^{2N-2m+2}(t_i-1)!\,(N+1-t_i)!}\\
\times
{\prodl _{1\le
i<j\le N-m+1} ^{}(t_{2j}-t_{2i})\prodl _{1\le
i<j\le N-m+1} ^{}(t_{2j-1}-t_{2i-1})}.
\endmultline\tag4.2$$
\endproclaim
\demo{Proof}Let $S=\{s_1,s_2,\dots,s_{2m-N-1}\}$ be the complement of
$T=\{t_1,t_2,\dots,\mathbreak t_{2N-2m+2}\}$ in $\{1,2,\dots,N+1\}$. 

Any perfect matching of the holey Aztec rectangle under
consideration can be naturally split along the central horizontal
line. The lower half is a perfect matching of an $m\times N$
Aztec rectangle, where $m$ vertices in the top-most row 
have been removed, including the vertices whose indices are in $S$.
Likewise, the upper half is a perfect matching of an $m\times N$
Aztec rectangle, where $m$ vertices in the bottom-most row 
have been removed, including the vertices whose indices are in $S$,
and avoiding the other vertices that are removed in the lower half.
Hence, by Lemma~2, the total number of perfect matchings that we are
interested in is the sum
$$\sum _{} ^{}\frac {2^{\binom {m}2}} {\prodl _{i=1} ^{m}(i-1)!}\prodl _{1\le
i<j\le m} ^{}(a_j-a_i)\frac {2^{\binom {m}2}} {\prodl _{i=1} ^{m}(i-1)!}\prodl _{1\le
i<j\le m} ^{}(b_j-b_i),$$
where the sum is over all pairs of sets $A=\{a_1,a_2,\dots,a_m\}$ and
$B=\{b_1,b_2,\dots,b_m\}$, whose union is
the complete set of indices $\{1,2,\dots,N+1\}$, and whose intersection
is $S$. Next we
extract the common factors that are generated by $A$ and $B$, thus
obtaining the sum 
$$\multline
2^{m^2-m}\frac {\prodl _{i=1} ^{2m-N-1}(s_i-1)!\,(N+1-s_i)!} {\prod
_{i=1} ^{m}(i-1)!^2}\\
\times
\sum _{} ^{}\prodl _{1\le
i<j\le N-m+1} ^{}(a'_j-a'_i)\prodl _{1\le
i<j\le N-m+1} ^{}(b'_j-b'_i),
\endmultline$$
where the sum is over all pairs of disjoint sets
$A'=\{a'_1,a'_2,\dots,a'_{N-m+1}\}$ and
$B=\{b'_1,b'_2,\dots,b'_{N-m+1}\}$, whose union is
$\{1,2,\dots,N+1\}\backslash S$.

Now we may use (3.1) to rewrite this last expression as
$$\frac {2^{m^2-m}\prodl _{i=1} ^{2m-N-1}(s_i-1)!\,(N+1-s_i)!} {\prod
_{i=N-m+2} ^{m}(i-1)!^2}
\sum _{} ^{}s_{\la(A')}(1^{N-m+1})\cdot s_{\la(B')}(1^{N-m+1}).
$$
The sum can be evaluated by means of Theorem~3. Renewed use of (3.1)
and some simplification eventually leads to (4.2).\quad \quad
\qed
\enddemo

Our next two theorems concern the enumeration of perfect matchings of
holey Aztec rectangles, where the ``holes" are on a horizontal row
next to the central horizontal row. 
Both theorems, when suitably specialized, solve Propp's Problem~10 in
\cite{\PropAA} that was mentioned in the Introduction. 
The second theorem, Theorem~10, 
was also obtained independently by Helfgott and Gessel
\cite{\HeGeAA, Theorem~1, second part}.

In general, we say that a horizontal row, $H$ say, is {\it by $d$ below the
central horizontal row\/}, if it is below the central horizontal row
and if a shortest path from any vertex of $H$
to some vertex on the central horizontal row
along the edges of the
Aztec rectangle needs exactly $d$ steps.
\proclaim{Theorem~9}Let $m$ and $N$ be positive integers with $2m+1\le N$. Then 
the number of perfect matchings of a $(2m+1)\times N$
Aztec rectangle, where all the vertices on the horizontal row that
is by 1 below the central row, except for 
the $t_1$-st, the $t_2$-nd, \dots, and the $t_{2m+1}$-st vertex, have been
removed (see Figure~5.a for an example with $m=2$, $N=7$, $a_1=1$,
$a_2=2$, $a_3=4$, $a_4=5$, $a_5=7$), equals
$$\frac {2^{m^2+3m+1}} {\prodl _{i=1} ^{m}(i-1)!\prodl _{i=1} ^{m+1}(i-1)!}
\prodl _{1\le
i<j\le m} ^{}(t_{2j}-t_{2i})\prodl _{1\le
i<j\le m+1} ^{}(t_{2j-1}-t_{2i-1}).\tag4.3$$
\endproclaim
\demo{Proof} Let $T=\{t_1,t_2,\dots,t_{2m+1}\}$.

Similar to the proof of Theorem~8, we may split any perfect matching of
the holey Aztec rectangle under consideration along this horizontal
line which is by 1 off the central line. An application of Lemma~1
yields that the number of perfect matchings that we are interested in
equals
$$\sum _{} ^{}\frac {2^{\binom {m+1}2}} {\prodl _{i=1} ^{m}(i-1)!}\prodl _{1\le
i<j\le m} ^{}(a_j-a_i)\frac {2^{\binom {m+2}2}} {\prodl _{i=1} ^{m+1}(i-1)!}\prodl _{1\le
i<j\le m+1} ^{}(b_j-b_i),$$
where the sum is over all pairs of disjoint sets $A=\{a_1,a_2,\dots,a_m\}$ and
$B=\{b_1,b_2,\dots,b_{m+1}\}$, whose union is the set $T$.
By using (3.1), it is seen that this expression can be rewritten as
$$ {2^{m^2+2m+1}} 
\sum _{} ^{}s_{\la(A)}(1^{m})\cdot s_{\la(B)}(1^{m+1}),
$$
where the sum is over the same set of pairs $(A,B)$. Clearly, the sum
can be evaluated by means of Theorem~4. Renewed use of (3.1)
gives (4.3) immediately.\quad \quad \qed
\enddemo
\vskip10pt
\vbox{
$$
\Einheit.4cm
\PfadDicke{.5pt}
\thinlines
\Pfad(0,4),43434343434343\endPfad
\Pfad(0,4),34343434343434\endPfad
\Pfad(0,2),43434343434343\endPfad
\Pfad(0,2),34343434343434\endPfad
\Pfad(0,0),4343\endPfad
\Pfad(6,0),4343\endPfad
\Pfad(12,0),43\endPfad
\Pfad(0,0),34343434343434\endPfad
\Pfad(0,-2),43434343434343\endPfad
\Pfad(0,-2),3434\endPfad
\Pfad(6,-2),3434\endPfad
\Pfad(12,-2),34\endPfad
\Pfad(0,-4),43434343434343\endPfad
\Pfad(0,-4),34343434343434\endPfad
\DickPunkt(1,-1)
\DickPunkt(3,-1)
\Kreis(5,-1)
\DickPunkt(7,-1)
\DickPunkt(9,-1)
\Kreis(11,-1)
\DickPunkt(13,-1)
\DickPunkt(0,0)
\DickPunkt(2,0)
\DickPunkt(4,0)
\DickPunkt(6,0)
\DickPunkt(8,0)
\DickPunkt(10,0)
\DickPunkt(12,0)
\DickPunkt(14,0)
\DickPunkt(1,1)
\DickPunkt(3,1)
\DickPunkt(5,1)
\DickPunkt(7,1)
\DickPunkt(9,1)
\DickPunkt(11,1)
\DickPunkt(13,1)
\DickPunkt(0,2)
\DickPunkt(2,2)
\DickPunkt(4,2)
\DickPunkt(6,2)
\DickPunkt(8,2)
\DickPunkt(10,2)
\DickPunkt(12,2)
\DickPunkt(14,2)
\DickPunkt(1,3)
\DickPunkt(3,3)
\DickPunkt(5,3)
\DickPunkt(7,3)
\DickPunkt(9,3)
\DickPunkt(11,3)
\DickPunkt(13,3)
\DickPunkt(0,4)
\DickPunkt(2,4)
\DickPunkt(4,4)
\DickPunkt(6,4)
\DickPunkt(8,4)
\DickPunkt(10,4)
\DickPunkt(12,4)
\DickPunkt(14,4)
\DickPunkt(1,5)
\DickPunkt(3,5)
\DickPunkt(5,5)
\DickPunkt(7,5)
\DickPunkt(9,5)
\DickPunkt(11,5)
\DickPunkt(13,5)
\DickPunkt(1,-3)
\DickPunkt(3,-3)
\DickPunkt(5,-3)
\DickPunkt(7,-3)
\DickPunkt(9,-3)
\DickPunkt(11,-3)
\DickPunkt(13,-3)
\DickPunkt(0,-2)
\DickPunkt(2,-2)
\DickPunkt(4,-2)
\DickPunkt(6,-2)
\DickPunkt(8,-2)
\DickPunkt(10,-2)
\DickPunkt(12,-2)
\DickPunkt(14,-2)
\DickPunkt(1,-5)
\DickPunkt(3,-5)
\DickPunkt(5,-5)
\DickPunkt(7,-5)
\DickPunkt(9,-5)
\DickPunkt(11,-5)
\DickPunkt(13,-5)
\DickPunkt(0,-4)
\DickPunkt(2,-4)
\DickPunkt(4,-4)
\DickPunkt(6,-4)
\DickPunkt(8,-4)
\DickPunkt(10,-4)
\DickPunkt(12,-4)
\DickPunkt(14,-4)
\hbox{\hskip8cm}
\Pfad(0,5),434343\endPfad
\Pfad(0,5),343434\endPfad
\Pfad(0,3),434343\endPfad
\Pfad(0,3),343434\endPfad
\Pfad(0,1),434343\endPfad
\Pfad(0,1),343434\endPfad
\Pfad(1,0),43\endPfad
\Pfad(1,-2),34\endPfad
\Pfad(0,-3),434343\endPfad
\Pfad(0,-3),343434\endPfad
\Pfad(0,-5),434343\endPfad
\Pfad(0,-5),343434\endPfad
\DickPunkt(1,-2)
\DickPunkt(3,-2)
\DickPunkt(5,-2)
\Kreis(0,-1)
\DickPunkt(2,-1)
\Kreis(4,-1)
\Kreis(6,-1)
\DickPunkt(1,0)
\DickPunkt(3,0)
\DickPunkt(5,0)
\DickPunkt(0,1)
\DickPunkt(2,1)
\DickPunkt(4,1)
\DickPunkt(6,1)
\DickPunkt(1,2)
\DickPunkt(3,2)
\DickPunkt(5,2)
\DickPunkt(0,3)
\DickPunkt(2,3)
\DickPunkt(4,3)
\DickPunkt(6,3)
\DickPunkt(1,4)
\DickPunkt(3,4)
\DickPunkt(5,4)
\DickPunkt(0,5)
\DickPunkt(2,5)
\DickPunkt(4,5)
\DickPunkt(6,5)
\DickPunkt(1,6)
\DickPunkt(3,6)
\DickPunkt(5,6)
\DickPunkt(1,-4)
\DickPunkt(3,-4)
\DickPunkt(5,-4)
\DickPunkt(0,-3)
\DickPunkt(2,-3)
\DickPunkt(4,-3)
\DickPunkt(6,-3)
\DickPunkt(1,-6)
\DickPunkt(3,-6)
\DickPunkt(5,-6)
\DickPunkt(0,-5)
\DickPunkt(2,-5)
\DickPunkt(4,-5)
\DickPunkt(6,-5)
\hskip3.5cm
$$
\line{\eightpoint\hskip.1cm a. A $5\times 7$ Aztec rectangle with some 
vertices\hskip.5cm
b. A $6\times 3$ Aztec rectangle with some vertices\hss}
\line{\eightpoint\hphantom{\hskip.1cm a. }on the row next to the central row removed\hskip.8cm
on the row next to the central line removed\hss}
\vskip10pt
\centerline{\eightpoint Figure 5}
}
\vskip10pt
\proclaim{Theorem~10}Let $m$ and $N$ be positive integers with $2m\ge N$. Then 
the number of perfect matchings of a $2m\times N$
Aztec rectangle, where all the vertices on the horizontal row that
is by 1 below the central row, except for 
the $t_1$-st, the $t_2$-nd, \dots, and the $t_{2N-2m+1}$-st vertex,
have been removed (see Figure~5.b for an example with $m=3$, $N=3$,
$a_1=2$), equals
$$\multline
{2^{m^2-m+N}\prodl _{i=m+1} ^{N+1}(i-1)!\prodl _{i=m+2} ^{N+1}(i-1)!}\\
\times
\frac {\prodl _{1\le
i<j\le N-m} ^{}(t_{2j}-t_{2i})\prodl _{1\le
i<j\le N-m+1} ^{}(t_{2j-1}-t_{2i-1})} {\prodl _{i=1}
^{2N-2m+1}(t_i-1)!\,(N+1-t_i)!}.
\endmultline\tag4.4$$
\endproclaim
Theorem~10 can be proved in the same way as Theorem~8, except that
Theorem~4 is used instead of Theorem~3. We leave the
details to the reader. 

Of course,
the special case $N=2m+2$ of Theorem~9 as well as the special case
$N=2m-1$ of Theorem~10 yields a solution of Problem~10 in
\cite{\PropAA}.

Clearly, we may continue in this manner, and derive formulas for 
the number of perfect matchings of an $M\times N$
Aztec rectangle, where $\vert N-M\vert$ vertices on a horizontal
row which is by $d$ off the central horizontal row
have been removed. By following the line of the proofs of
Theorems~8--10,
we obtain $\fl{d/2}$-fold summations if we apply Theorem~6. We just
state the result and leave the routine details of verification to
the reader.
\proclaim{Theorem~11}Let $m$ and $N$ be positive integers, and let $d$ be
a nonnegative integer, with $2m+d\le N$. Then 
the number of perfect matchings of a $(2m+d)\times N$
Aztec rectangle, where all the vertices on the horizontal row that
is by $d$ below the central row, except for 
the $t_1$-st, the $t_2$-nd, \dots, and the $t_{2m+d}$-th vertex, have been
removed, equals
$$\multline
\frac {2^{m^2+(d+2)m+\binom {d+1}2}} {\prodl _{i=1} ^{m}(i-1)!\prodl _{i=1}
^{m+d}(i-1)!}
\prodl _{1\le
i<j\le m+\fl{d/2}} ^{}(t_{2j}-t_{2i})
\prodl _{1\le
i<j\le m+\cl{d/2}} ^{}(t_{2j-1}-t_{2i-1})\\
\times
\sum _{1\le k_1<k_2<\dots<k_{\fl{d/2}}\le m+\fl{d/2}} ^{}\Bigg(
\prodl _{1\le i<j\le \fl{d/2}} ^{}(t_{2k_j}-t_{2k_i})^2
\prodl _{i=1} ^{\fl{d/2}}\frac {\prodl _{j=1} ^{m+\cl{d/2}}
\v{t_{2k_i}-t_{2j-1}}}
{\underset j\ne k_i\to{\prodl _{j=1}
^{m+\fl{d/2}}}\v{t_{2k_i}-t_{2j}}}\Bigg).\\
\endmultline\tag4.5$$
\endproclaim
\proclaim{Theorem~12}Let $m$ and $N$ be positive integers, and let $d$ be
a nonnegative integer, with $2m+d-1\ge N$. Then 
the number of perfect matchings of a $(2m+d-1)\times N$
Aztec rectangle, where all the vertices on the horizontal row that
is by $d$ below the central row, except for 
the $t_1$-st, the $t_2$-nd, \dots, and the $t_{2N-2m-d+2}$-nd vertex, have been
removed, equals
$$\multline
{2^{m^2+(d-2)m+\binom {d-1}2+N}}\frac {\prodl _{i=m+1} ^{N+1}(i-1)!\prodl
_{i=m+d+1} ^{N+1}(i-1)!} 
{\prodl _{i=1}
^{2N-2m-d+2}(t_i-1)!\,(N+1-t_i)!}\\
\times
\prodl _{1\le
i<j\le N-m+1-\cl{d/2}} ^{}(t_{2j}-t_{2i})
\prodl _{1\le
i<j\le N-m+1-\fl{d/2}} ^{}(t_{2j-1}-t_{2i-1})\\
\times
\sum _{1\le k_1<k_2<\dots<k_{\fl{d/2}}\le N-m+1-\cl{d/2}} ^{}\Bigg(
\prodl _{1\le i<j\le \fl{d/2}} ^{}(t_{2k_j}-t_{2k_i})^2\hskip3cm\\
\cdot
\prodl _{i=1} ^{\fl{d/2}}\frac {\prodl _{j=1} ^{N-m+1-\fl{d/2}}
\v{t_{2k_i}-t_{2j-1}}}
{\underset j\ne k_i\to{\prodl _{j=1}
^{N-m+1-\cl{d/2}}}\v{t_{2k_i}-t_{2j}}}\Bigg).
\endmultline\tag4.6$$
\endproclaim
It is natural to ask if there are cases where the sums in (4.5) or
(4.6) can be simplified. Indeed, thanks to Theorem~6, there are two such 
cases. One is the case where the gaps between the vertices that are not removed 
are always the same,
i.e., where the corresponding $t_i$'s form an arithmetic progression.
The other is the case where the quotients of successive gaps 
between the vertices that are not removed 
are always the same,
i.e., where the corresponding $t_i$'s form a (shifted) geometric progression.

The following two theorems contain the results for the case of
arithmetic progressions. We want
to point out that both of them (as well as Theorems~7 and 8 of course)
contain the well-known enumeration of all perfect matchings of an Aztec
diamond \cite{\ElKLAA} as special case.
\proclaim{Theorem~13}Let $m,N,C,D$ be positive integers, and let $d$
be a nonnegative integer, with $2m+d\le N$ and
$C+(2m+d-1)D\le N$. Then 
the number of perfect matchings of a $(2m+d)\times N$
Aztec rectangle, where all the vertices on the horizontal row that
is by $d$ below the central row, except for 
the $C$-th, the $(C+D)$-th, $(C+2D)$-th, \dots, and the 
$(C+(2m+d-1)D)$-th vertex, have been removed, equals
$$
2^{\binom {2m+d+1}2}\,D^{m^2+(d-1)m+\binom d2}.
\tag4.7$$
\endproclaim
\proclaim{Theorem~14}Let $m,N,C,D$ be positive integers, and let $d$
be a nonnegative integer, with $2m+d-1\le N$ and
$C+(2N-2m-d+1)D\le N+1$. Then 
the number of perfect matchings of a $(2m+d-1)\times N$
Aztec rectangle, where all the vertices on the horizontal row that
is by $d$ below the central row, except for 
the $C$-th, the $(C+D)$-th, $(C+2D)$-th, \dots, and the 
$(C+(2N-2m-d+1)D)$-th vertex, have been removed, equals
$$\multline
2^{\binom {2m+d}2+(N+1)(N-2m-d+1)}\,D^{m^2+(d-1)m+\binom
d2+N(N-2m-d+1)}\\
\times
\frac {\prodl _{i=m+1} ^{N+1}(i-1)!\prodl _{i=m+d+1} ^{N+1}(i-1)!
\prodl _{i=1} ^{N-m+1}(i-1)!\prodl _{i=1} ^{N-m-d+1}(i-1)!} 
{\prodl _{i=1} ^{2N-2m-d+2}(C+Di-1)!\,(N+1-C-Di)!}.
\endmultline\tag4.8$$
\endproclaim
\demo{Proof of Theorems~13 and 14} We set $t_i=C+D(i-1)$, $i=1,2,\dots$,
in Theorems~11 and 12. The resulting sums turn out to be exactly of
the form of the left-hand side of (3.5), with $x=1/2$ and $y=3/2$,
respectively $x=3/2$ and $y=3/2$, depending on $d$ being even or odd. 
Hence, it can be evaluated.
The obtained expressions can be drastically simplified and eventually
turn into (4.7) and (4.8), respectively.\quad \quad \qed
\enddemo

The results for the case of (shifted) geometric progressions are the
following two.
\proclaim{Theorem~15}Let $m$ and $N$ be positive integers, let $d$ be
a nonnegative integer, let $C,D,q$ be rational numbers, $q>1$, such
that $(C+D)$, $C+Dq$, $C+Dq^2$, 
\dots, $C+Dq^{2m+d-1}$ are integers, and such that $2m+d\le N$ and
$C+Dq^{2m+d-1}\le N$. Then 
the number of perfect matchings of a $(2m+d)\times N$
Aztec rectangle, where all the vertices on the horizontal row that
is by $d$ below the central row, except for 
the $(C+D)$-th, the $(C+Dq)$-th, $(C+Dq^2)$-th, \dots, and the 
$(C+Dq^{2m+d-1})$-th vertex, have been removed, equals
$$\multline
2^{m^2+(d+2)m+\binom {d+1}2}\,D^{m^2+(d-1)m+\binom d2}
{q^{\frac {1} {6}(m+d-1)(4m^2+(2d+1)m+d(d-2))}}\\
\times
\frac {\prodl _{i=1} ^{m}\v{(q^2;q^2)_{i-1}}
\prodl _{i=1} ^{m+d}\v{(q^2;q^2)_{i-1}}} {
{\prodl _{i=1} ^{d}(-q;q)_{i-1}}\prodl _{i=1} ^{m}(i-1)!
\prodl _{i=1} ^{m+d}(i-1)!}.
\endmultline\tag4.9$$
\endproclaim
\proclaim{Theorem~16}Let $m$ and $N$ be positive integers, let $d$ be
a nonnegative integer, let $C,D,q$ be rational numbers, $q>1$, such
that $(C+D)$, $C+Dq$, $C+Dq^2$, 
\dots, $C+Dq^{2N-2m-d+1}$ are integers, and such that $2m+d-1\le
N$ and $C+Dq^{2N-2m-d+1}\le N+1$. Then 
the number of perfect matchings of a $(2m+d-1)\times N$
Aztec rectangle, where all the vertices on the horizontal row that
is by $d$ below the central row, except for 
the $(C+D)$-th, the $(C+Dq)$-th, $(C+Dq^2)$-th, \dots, and the 
$(C+Dq^{2N-2m-d+1})$-th vertex, have been removed, equals
$$\multline
{2^{m^2+(d-2)m+\binom {d-1}2+N}}\,D^{m^2+(d-1)m+\binom
d2+N(N-2m-d+1)}\\
\times
{q^{\frac {1} {6}(N-m)(4(N-m-d+1)^2+(2d+1)(N-m-d+1)+d(d-2))}}\\
\times
\frac {\prodl _{i=m+1} ^{N+1}(i-1)!\prodl
_{i=m+d+1} ^{N+1}(i-1)!
\prodl _{i=1} ^{N-m+1}\v{(q^2;q^2)_{i-1}}\prodl _{i=1} ^{N-m-d+1}\v{(q^2;q^2)_{i-1}}} 
{\prodl _{i=1} ^{d}(-q;q)_{i-1}\prodl _{i=1}
^{2N-2m-d+2}(C+Dq^{i-1}-1)!\,(N+1-C-Dq^{i-1})!}.
\endmultline\tag4.10$$
\endproclaim
\demo{Proof of Theorems~15 and 16} We set $t_i=C+Dq^{i-1}$, $i=1,2,\dots$,
in Theorems~11 and 12. The resulting sums turn out to be exactly of
the form of the left-hand side of (3.6), with $q\to q^2$, $x=q$ and
$y=q^3$,
respectively $x=q^3$ and $y=q^3$, depending on $d$ being even or odd. 
Hence, it can be evaluated.
The obtained expressions can be drastically simplified and eventually
turn into (4.9) and (4.10), respectively.\quad \quad \qed
\enddemo

It should be observed that, formally, Theorems~13 and 14 can be seen
as limit cases of Theorems~15 and 16. Namely, if in Theorems~15 and
16 we replace $C$ by $\frac {q^C-1} {q-1}-\frac {q^D-1} {(q-1)^2}$
and $D$ by $\frac {q^D-1} {(q-1)^2}$, and then perform the limit
$q\to 1^+$, we obtain, formally, Theorems~13 and 14.

\enddocument